\documentclass[12pt]{article}

\usepackage{amsthm,color}
\usepackage{amsfonts}
\usepackage{graphicx}
\usepackage{mathrsfs}
\usepackage{amsmath}
\usepackage{amssymb}
\usepackage{float}
\usepackage{natbib}
\usepackage{booktabs}
\usepackage{url}
\usepackage{multirow}
\usepackage{arydshln}
\usepackage{threeparttable}
\usepackage{courier}
\usepackage{authblk} 

\usepackage{multicol}
\usepackage{latexsym}
\usepackage{psfrag}
\usepackage[usenames,dvipsnames]{xcolor}

\usepackage{breakcites}
\usepackage{bbm}
\usepackage{epsfig,epstopdf}
\usepackage[colorlinks, linkcolor=black, citecolor=black]{hyperref}

\usepackage[utf8]{inputenc}
\usepackage[english]{babel}
\usepackage{caption}
\usepackage{setspace}

\usepackage{sectsty}
\sectionfont{\fontsize{12}{15}\selectfont}
\subsectionfont{\fontsize{12}{15}\selectfont}

\newtheorem{theorem}{Theorem}

\newtheorem{lemma}[theorem]{Lemma}
\newtheorem*{remark}{Remark}

\def\A{\boldsymbol{A}}
\def\B{\boldsymbol{B}}

\def\J{\boldsymbol{J}}
\def\C{\boldsymbol{C}}
\def\D{\boldsymbol{D}}
\def\q{\boldsymbol{q}}
\def\Q{\boldsymbol{Q}}

\def\g{\boldsymbol{g}}
\def\G{\boldsymbol{G}}

\def\e{\boldsymbol{e}}
\newcommand{\bu}{\boldsymbol u}
\def\m{\boldsymbol{m}}

\def\G{\mathcal{G}}
\newcommand{\diag}{\rm\mbox{diag}}

\def\bGamma{\boldsymbol{\Gamma}}
\def\bgamma{\boldsymbol{\gamma}}

\def\bbeta{\boldsymbol{\beta}}

\def\btheta{\boldsymbol{\theta}}
\def\bstheta{\mbox{\scriptsize \boldmath $\theta$}}
\def\htheta{\hat{\btheta}}
\def\etab{\boldsymbol{\eta}}

\def\hetab{\hat{\etab}}

\def\bSigma{\boldsymbol{\Sigma}}
\def\bnu{\boldsymbol{\nu}}
\def\bsnu{\mbox{\scriptsize \boldmath $\nu$}}

\def\Sumij{\sum_{i=0}^1 \sum_{j=1}^{n_{i1}}}

\def\ba{\left( \begin{array}}
\def\ea{\end{array} \right)}

\def\lb{\left\{ }
\def\rb {\right\} }
\def\Lm{\left[ }
\def\Rm{\right] }

\begin{document}

{\centering {\large {\bf Semiparametric inference on Gini indices  of two  semicontinuous populations under density ratio models}}\par }
\bigskip

\centerline{Meng Yuan, \ Pengfei Li \ and \ Changbao Wu\footnote{Meng Yuan is  doctoral student, Pengfei Li is Professor and Changbao Wu is Professor, Department of Statistics and Actuarial Science, University of Waterloo, Waterloo ON N2L 3G1,Canada (E-mails: {\em m33yuan@uwaterloo.ca}, \ {\em pengfei.li@uwaterloo.ca} \ and \ {\em cbwu@uwaterloo.ca}).}}

\bigskip

\bigskip

\hrule

{\small
\begin{quotation}
\noindent
The Gini index is a popular  inequality measure with many applications in social and economic studies. 
This paper studies semiparametric inference on the Gini indices of two semicontinuous populations. 
We characterize the distribution of each semicontinuous population by a mixture of a discrete point mass at zero and a continuous skewed positive component.
A semiparametric density ratio model is then employed to link the positive components of the two distributions. 
We propose the maximum empirical likelihood estimators of the two Gini indices and their difference,  
and further investigate  the asymptotic properties of the proposed estimators.
The asymptotic results enable us to construct confidence intervals and perform hypothesis tests
for the two Gini indices and their difference. 
We show that the proposed estimators
are more efficient than the existing fully nonparametric estimators. 
The proposed estimators and the asymptotic results are also applicable to cases without excessive zero values. 
Simulation studies show the superiority of our proposed method over existing methods.  
Two real-data applications are presented using the proposed methods.

\vspace{0.3cm}

\noindent
{\bf Keywords:}  Density ratio model, zero-excessive data, empirical likelihood, equality test.
\end{quotation}
}

\hrule

\bigskip

\bigskip

\section{Introduction}
\label{intro}

The Gini index, first proposed by \cite{gini1912variabilita},  has been widely used to measure population inequality. 
In economic studies, it is an important measure of income or wealth inequality among
individuals or households in a particular population \citep{wang2016jackknife,peng2011empirical}. In life expectancy studies, 
it is used to describe the concentration of survival times and to evaluate inequality among people in the target population \citep{bonetti2009gini, lv2017gini}.
The index is closely related to the Lorenz curve \citep{lorenz1905methods}, a widely used measure for the size distribution of income or wealth. 
It is the ratio of the area between the Lorenz curve and the 45-degree line to the area under the 45-degree line.  
Hence, the Gini index ranges from 0 to 1, with 0 indicating perfect equality and 1 for extreme inequality. 

Study variables such as income and survival time are often modelled by using a positive continuous distribution. One important scenario in applications is that there are two related populations, each containing a sizeable zero values for the study variable. The inferential problems can be on the Gini index for each population separately or the difference of the two Gini indices. The scenario is quite common in practice but efficient inferential procedures are not available in the existing literature.

In this paper, we propose new semiparametric inference procedures for the Gini indices of two semicontinuous populations.
Specifically, suppose that we have two independent samples from two related populations with values of the study variable $X$ generated from the following mixture models:
\begin{equation}
\label{mixture}
  (X_{i1},\cdots, X_{in_i}) \sim F_i(x)= \nu_iI(x \geq 0)+(1-\nu_i)I(x>0)G_i(x), ~~~\mbox{for~} ~i = 0,1,
\end{equation}
where $\nu_i$ is the zero proportion in population $i$, $n_i$ is the sample size for sample $i$, $I(\cdot)$ is an indicator function, and $G_i(\cdot)$ is the cumulative distribution function (CDF) of the positive observations in sample $i$. 
For population $i =0,1$, the Gini index can be equivalently defined \citep{david1968miscellanea} as
\begin{equation}
\label{gini}
  \G_i = \frac{D_i}{2\mu_i},
\end{equation}
where 
$
D_i = E( |X_{i1}-X_{i2}|)
$
is the expected absolute difference of $X$ for two randomly selected units from population $i$ 
and 
$ \mu_i = E(X_{i1})
$
is the expectation  of population $i$.  Our discussions in this paper focus on statistical inferences on $ \G_0$, $ \G_1$, and $ \G_0- \G_1$. 
It is worth mentioning that although our results are presented for cases where the two populations contain excess zeros for the study variable, 
the proposed methods and the theoretical results are also applicable to cases without excess zeros, i.e., $\nu_i=0$ in model \eqref{mixture}. In addition, inferences on a general function of $ \G_0$ and $ \G_1$ can also be conducted. See Section~\ref{estimation} for further discussion. 

Samples with positive outcomes only, i.e., $\nu_i=0$ in model \eqref{mixture}, 
are common in studies of family income or wealth or a country's gross domestic products \citep{gastwirth1972estimation,cowell2011measuring}. 
For instance, all the household incomes are positive in the 1997 Family and Income and Expenditure Survey 
conducted by the Philippine Statistics Authority. More details can be found in Section~\ref{realdata}. 
Samples with a mixture of excess zero values and skewed positive outcomes,  i.e., $\nu_i>0$ in model \eqref{mixture}, 
naturally arise in studies of expenditure data and health cost data \citep{ZhouTu1999, zhou2000interval, zhou2008computer}. 
For example, \cite{zhou2008computer} presented a dataset from the assessment of inpatient charges (see Section~\ref{realdata}), and 
most patients with uncomplicated hypertension had no hospitalization and therefore zero costs.
This paper systematically studies both cases ($\nu_i=0$ and $\nu_i>0$) in a unified framework via model \eqref{mixture}. 

Many studies of the Gini index have applied nonparametric methods. 
For example, point estimators of $ \G_0$, $ \G_1$, and $ \G_0- \G_1$ and their asymptotic variance estimation 
have been discussed in \cite{Hoeffding1948asymptotic},  \cite{anand1983inequality}, 
\cite{ogwang2000convenient},  \cite{giles2004calculating}, 
\cite{modarres2006cautionary},  
\cite{yitzhaki1991calculating}, \cite{karagiannis2000method},  
and 
\cite{davidson2009reliable}. 
See \cite{wang2016jackcompare} for a detailed review. 
\cite{qin2010empirical} and \cite{peng2011empirical}   used  the empirical likelihood (EL) method \citep{owen} to 
construct confidence intervals (CIs) for the index. 
More recently, \cite{wang2016jackknife} derived the jackknife EL (JEL). 
\cite{peng2011empirical} and \cite{wang2016jackcompare} compared two Gini indices 
of independent or correlated populations using the EL method and the JEL method, respectively. 


Fully nonparametric methods are robust to potential model misspecifications on the $F_i$'s. However, these methods ignore the characteristics common to the two samples and/or the relation between the two populations, which have been shown to be useful for more efficient statistical inferences; see, for instance, the  studies on the strengths of lumber produced in Canada in different years
\citep{chen2013quantile,cai2017hypothesis,cai2018empirical}  and biomarkers for the diagnosis of disease in case and control groups \citep{yuan2020zeros, yuan2021semiparametric}.

To combine the information from the two samples without making risky parametric distributional assumptions, 
we model the CDFs of the positive observations $G_0$ and $G_1$ in model \eqref{mixture} via a density ratio model (DRM) \citep{Anderson1979,Qin2017}. 
Let $dG_i$ be  the probability density function of $G_i$, $i=0,1$.
The DRM postulates that 
\begin{equation}
\label{drm}
dG_1(x) = \exp\{\alpha +\boldsymbol{\beta}^\top\q(x)\}dG_0(x) = \exp\{\btheta^\top\Q(x)\}dG_0(x)\,,
\end{equation}
where $\btheta = (\alpha,\bbeta^\top)^\top$ is the vector of unknown parameters and $\Q(x) = (1,\q(x)^\top)^\top$ with $\q(x)$ being a $d$-dimensional, prespecified,  nontrivial basis function.
Note that the baseline distribution $G_0$ is completely unspecified, which allows the DRM to embrace many distributions that are commonly used in the study of the Gini index \citep{moothathu1985distributions,moothathu1989unbiased,moothathu1990best}.  
For example, when $\q(x) = \log(x)$, the DRM includes two log-normal distributions with the same variance with respect to the log-scale,
two Pareto distributions with the same scale parameter, 
and two chi-square distributions with different degrees of freedom;
when $\q(x) = x$, the DRM includes 
two exponential distributions with different rates. 
The DRM is closely related to the well-studied logistic regression \citep{qin1997goodness} 
and is broader than the Cox proportional hazard model \citep{jiang2012inference}.

The DRM  has served as a useful and flexible inferential platform for many statistical problems,
including quantile and quantile-function estimation \citep{chen2013quantile,yuan2021semiparametric,chen2021composite},
receiver operating characteristic analysis \citep{qin2003using,chen2016using,yuan2021semiparametric}, hypothesis-testing problems \citep{fokianos2001semiparametric,cai2017hypothesis}, and data integration \citep{qin2014using,yuan2021semiparametric}.
Recently,  the DRM has been used for semicontinuous data. 
For example, 
\cite{Wang2017} and \cite{Wang2018} developed EL ratio (ELR) statistics for testing the homogeneity of distributions and the equality of population means, respectively. \cite{yuan2020zeros} proposed estimators of linear functionals of two semicontinuous populations.  
To the best of our knowledge, inferential procedures for two Gini indices $\G_0$, $\G_1$ and their difference $\G_0 - \G_1$ 
have not been explored under the mixture model \eqref{mixture} and the DRM \eqref{drm}. 
This paper aims to fill this void.

Our contributions can be summarized as follows. 
We first propose the maximum EL estimators (MELEs) of  $\G_0$, $\G_1$, and $\G_0 - \G_1$  under models \eqref{mixture} and \eqref{drm}. 
Using techniques from U-statistics and   V-statistics, 
we derive the asymptotic normality of the MELEs of   $\G_0$, $\G_1$, and $\G_0 - \G_1$
and show that their asymptotic variances are smaller than those of nonparametric estimators. 
These asymptotic results are used to construct Wald-type CIs for $\G_0$, $\G_1$, and $\G_0 - \G_1$
and obtain a Wald-type statistic for testing $H_0:\G_0 =\G_1$. 
The proposed methods and asymptotic results
are applicable to both $\nu_i=0$ and $\nu_i>0$ in model \eqref{mixture}.
Extensive simulation studies and applications to two real datasets demonstrate the advantages of our proposed methods over existing ones. 
Software for implementing the proposed methods  and \cite{qin1997goodness}'s goodness-of-fit test for the DRM assumption in  \eqref{drm}
has been developed in R and is available from the authors upon request.

%
%

The rest of the paper is organized as follows. In Section~\ref{estimation}, 
we propose the MELEs of $\G_0$, $\G_1$, and $\G_0 - \G_1$ and study their asymptotic properties under the assumed semiparametric density ratio model. 
We construct CIs and conduct hypothesis tests on $\G_0$, $\G_1$, and $\G_0 - \G_1$ based on the theoretical results. 
Results from simulation studies are presented in Section~\ref{simu},
and applications to two real-world datasets are given in  Section~\ref{realdata}. 
We conclude the paper with a discussion in Section~\ref{conclude}. 
Proofs and technical details and additional simulation results are provided in the Supplementary Material.

\section{Main Results} 
\label{estimation}

Let $n_{i0}$ and $n_{i1}$ be the (random) numbers of zero observations and positive observations, respectively, in each sample $i=0,1$.  
It follows that $n_i = n_{i0}+n_{i1}$ for $i =0,1$.
Without loss of generality, we assume that the first $n_{i1}$ observations in group $i$, $X_{i1},\cdots, X_{in_{i1}}$, are positive,
and the remaining $n_{i0}$ observations are 0.
Let $n$ be the total (fixed) sample size, i.e., $n = n_0+n_1$.

\subsection{The MELEs of  $\G_0$ and $\G_1$} 

In this section, we  develop the EL function and present the MELEs  of  $\G_0$ and $\G_1$.  
Based on the two samples from model \eqref{mixture}, 
the full likelihood function  is given by
\begin{eqnarray}
\label{gini.like}
L_n &=& \prod_{i=0}^1 \Lm \lb v_i^{n_{i0}}\left(1-v_i\right)^{n_{i1}} \rb \prod_{j=1}^{n_{i1}}dG_i\left(X_{ij}\right)\Rm,
\end{eqnarray}
which is the product of two likelihood components: one from the number of zero observations and the other from the positive observations. 

Following the EL principle \citep{owen} and with the help of the DRM \eqref{drm}, 
we use the combined sample to  model $G_0$ via 
\begin{equation} \label{gini.G0}
G_0(x)= \sum_{i=0}^1\sum_{j = 1}^{n_{i1}}p_{ij}I(X_{ij}\le x),
\end{equation}
where $p_{ij} = dG_0(X_{ij})$ for $i=0,1$ and $j=1,\ldots,n_{i1}$.
The model for $G_0$ in \eqref{gini.G0} and the DRM \eqref{drm} together imply that 
\begin{equation} \label{gini.G1}
G_1(x)= \sum_{i=0}^1\sum_{j = 1}^{n_{i1}}p_{ij} \exp\{\btheta^\top\Q(X_{ij})\} I(X_{ij}\le x).
\end{equation}
To ensure that both $G_0$ and $G_1$ are CDFs, 
the feasible  $p_{ij}$'s must satisfy the following set of constraints: 
\begin{equation}
\label{gini.cons}
\mathcal{C} = \left\{ (\btheta,\boldsymbol{P}): p_{ij} > 0,~
\sum_{i=0}^1\sum_{j = 1}^{n_{i1}}p_{ij}=1,~
\sum_{i=0}^1\sum_{j = 1}^{n_{i1}} p_{ij}\exp\{\btheta^\top\Q(X_{ij})\}=1 \right\}, 
\end{equation}
where $\boldsymbol{P} = \{p_{ij}: j=1,\ldots,n_{i1}, i=0,1\}$. 

Let $\bnu = (\nu_0,\nu_1)^\top$. Substituting \eqref{gini.G0} and \eqref{gini.G1} into \eqref{gini.like} and taking the logarithm, 
we obtain the empirical log-likelihood function of $(\boldsymbol{\nu},\boldsymbol{\theta}, \boldsymbol{P})$ 
as 
\begin{eqnarray}
\label{gini.full}
\tilde \ell(\bnu,\btheta, \boldsymbol{P} ) = \ell_0\left(\bnu\right) +  \tilde \ell_1\left(\boldsymbol{\theta},\boldsymbol{P}\right),
\end{eqnarray}
where 
\begin{equation*}
    \ell_0\left(\bnu\right) = \sum_{i =0}^1\log\left\{ v_i^{n_{i0}}\left( 1-v_i\right)^{n_{i1}} \right\}  \;\;\;
    ~\text{and}~ \;\;\; 
    \tilde \ell_1\left(\btheta,\boldsymbol{P}\right)  = \Sumij\log p_{ij}+ \sum_{j=1}^{n_{11}}\left\{\btheta^\top\Q(X_{1j})\right\}.
\end{equation*}
Here $\ell_0\left(\bnu\right)$ is the binomial log-likelihood function corresponding to the zero observations, and $\tilde \ell_1\left(\boldsymbol \theta,\boldsymbol{P}\right)$ represents the empirical log-likelihood function associated with the positive observations. 
The MELE of $(\boldsymbol{\nu},\boldsymbol{\theta}, \boldsymbol{P})$ is then defined as
\begin{equation*}
    (\hat{\bnu},\hat{\btheta}, \hat{\boldsymbol{P}}) = \arg\max_{\bnu,\btheta, \boldsymbol{P} } \tilde \ell(\bnu,\btheta, \boldsymbol{P} ) 
\end{equation*}
subject to the constraints in $\mathcal{C}$.

With the form of $\tilde \ell$ given in \eqref{gini.full},  it can be checked that  the MELE of $\bnu$ has a closed form expression given by 
$$
\hat{\bnu} = \arg\max_{\bnu} \ell_0(\bnu)=(n_{00}/n_0,n_{10}/n_1)^\top \,,
$$ 
and the MELE of $(\boldsymbol{\theta}, \boldsymbol{P})$ is then denoted as 
\begin{equation}
\label{max_theta_p}
(\hat \btheta,\hat{\boldsymbol{P}}) = \arg\max_{(\btheta, \boldsymbol{P})\in \mathcal{C}} \tilde \ell_1(\btheta, \boldsymbol{P} ) .
\end{equation} 

Following  \cite{cai2017hypothesis}, the estimator $\hat \btheta$ can be obtained by maximizing 
the following dual empirical log-likelihood function without any constraints: 
\begin{equation}
\label{gini.dual}
\ell_1(\btheta)=
- \Sumij
  \log\left\{ 1 + \hat{\rho}[\exp\{ \btheta^{\top} \Q(X_{ij})\}-1] \right\}
+ \sum_{j=1}^{n_{11}} \{\btheta^{\top}  \Q(X_{1j}) \},
\end{equation}
where $\hat{\rho} = n_{11}(n_{01}+n_{11})^{-1}$ is a random variable.
That is, $\hat{\btheta}= \arg \max_{\btheta}\ell_1(\btheta)$. Once the $\hat{\btheta}$ is obtained, 
the MELEs of the $ p_{ij} $'s are computed as  
\begin{equation}
\label{gini.pij}
    \hat p_{ij} = \frac{1}{n_{01}+n_{11}} \lb 1 + \hat{\rho}[\exp\{ \hat\btheta^{\top} \Q(X_{ij})\}-1]\rb ^{-1}.
\end{equation}
With the MELEs $\hat{\btheta}$ and $ \hat p_{ij}$,  
the MELE of $G_0(x)$ and $G_1(x)$ for $x >0$ as specified in \eqref{gini.G0}--\eqref{gini.G1} can be computed as 
\begin{equation}
\label{gini.hatG}
    \hat G_0(x) = \Sumij \hat p_{ij} I(X_{ij}\leq x)~~
    \mbox{ and }
    ~~
    \hat G_1(x) = \Sumij \hat p_{ij} \exp\{\hat \btheta^\top\Q(x)\} I(X_{ij}\leq x). 
\end{equation}

We now move to the point estimation of the Gini indices $\G_0$ and $\G_1$. 
In Section 1.1 of the Supplementary Material, we show that $\G_i$ in \eqref{gini}
can be equivalently expressed as 
\begin{equation}
\label{gini2}
\G_i=(2\nu_i-1)+(1-\nu_i){\psi_i}/{m_i},
\end{equation}
where $m_i = \int_0^{\infty} x dG_i(x) $ and $\psi_i = \int_0^{\infty} \{2xG_i(x)\}dG_i(x)$. 
Using the alternative form given in \eqref{gini2}, the MELE of the two Gini indices are given by 
\begin{eqnarray}
\label{hatGini}
\hat\G_i =(2\hat\nu_i-1)+(1-\hat\nu_i){\hat\psi_i}/{\hat m_i}, ~~i=0,1,
\end{eqnarray}
where 
\begin{eqnarray*}
\hat m_i=  \int_0^{\infty} x d \hat G_i(x) ~~
\mbox{and}
~~
\hat\psi_i=  \int_0^{\infty} \big\{2x \hat G_i(x)\big\}d \hat G_i(x). 
\end{eqnarray*}

\begin{remark}
We comment that the  MELEs of the two Gini indices  in \eqref{hatGini}
are also applicable to the case where  there is no excess of zero values, i.e., $\bnu = (0,0)$ and $n_{i1} = n_i$. 
We need to set $\hat\nu_i=0$ and obtain $\btheta$ by maximizing $\ell_1(\btheta)$ in \eqref{gini.dual} with $\hat{\rho} = n_{1}(n_{0}+n_{1})^{-1}$; 
then the MELEs in \eqref{gini.pij}--\eqref{hatGini} can be directly applied.  
%
\end{remark}

\subsection{Asymptotic properties of MELEs}

In this section, we study the asymptotic properties of the MELEs $(\hat\G_0,\hat\G_1)$. 
We use $\bnu^{*}$, $\btheta^{*}$, and  $(\G_0^*, \G_1^*)$  to denote the true values of $\bnu$, $\btheta$, and $(\G_0, \G_1)$, respectively. 
Let $w_i = n_i/n$ and  
\begin{gather*}
\Delta^*= \sum_{i=0}^1w_i(1-\nu_i^*),
~\rho^*=\frac{w_1(1-\nu_1^*)}{\Delta^*},~\omega(x)=\exp\{\btheta^{*\top}\Q(x)\},\\
h(x)=1+\rho^* \{\omega(x)-1\}, ~~~~h_1(x)=\rho^* \omega(x)/h(x),~~~~~~~~~~~~~~~~~~~\\
u_0(x) = (2\nu_0^*-1)x + (1-\nu_0^*)\left[2 \left\{xG_0(x)+ \int_{x}^{\infty} ydG_0(y)\right\} - \psi_0\right],\\
u_1(x) = (2\nu_1^*-1)x + (1-\nu_1^*)\left[2 \left\{xG_1(x)+ \int_{x}^{\infty} ydG_1(y)\right\} - \psi_1\right],\\
 \J =\ba{cccc}
    -\frac{\G_0^*}{m_0}&\frac{1}{m_0}&0&0\\
    0&0&-\frac{\G_1^*}{m_1}&\frac{1}{m_1}
    \ea,~
    \A_{\bstheta} =
\Delta^*(1-\rho^*)E_0\left\{h_1(X) \Q(X)\Q(X)^\top \right\},
\end{gather*}
where $E_0(\cdot)$ denotes the expectation operator with respect to $G_0$, and $X$ refers to a random variable from the distribution $G_0$. 

The asymptotic results in this section are developed under the following regularity conditions. 
\begin{itemize}
    \item [C1:] As the total sample size $n$ goes to infinity, $n_0/n \rightarrow w_0$  for some constant $w_0 \in (0,1)$.
    \item [C2:] The two CDFs $G_0$ and $G_1$ satisfy the DRM \eqref{drm} with the true parameter $\btheta^*$, and $\int_0^{\infty} \exp\{\btheta^\top\Q(x)\}dG_0(x) < \infty$ for all $\btheta$ in a neighborhood of the true value $\btheta^*$.
    \item [C3:] The components of $\Q(x)$ are continuous and stochastically linearly independent.
    \item [C4:] The moments $\int_0^{\infty} x^2dG_0(x)$ and $\int_0^{\infty}x^2\exp\{\btheta^\top\Q(x)\}dG_0(x)$ exist for all $\btheta$ in a neighborhood of the true value $\btheta^*$. 
\end{itemize}

Condition C1 indicates that both $n_0$ and $n_1$ go to infinity at the same rate.
For simplicity, and convenience of presentation, 
we write $w_0=n_0/n$ and assume that it is a constant. 
This does not affect our technical development.
Condition C2 guarantees the existence of the moment generating function of $\Q(X)$ in a neighborhood of $\btheta^*$  and therefore all its finite moments.
Condition C3 is an identifiability condition. 
Conditions C2 and C3 together imply that $\A_{\bstheta}$ is positive definite and the quadratic approximation of the dual empirical log-likelihood function $\ell_1(\btheta)$ is applicable.
Conditions C1--C4 guarantee that the linear approximations of  $\hat\G_0$ and $\hat\G_1$ can be used.

The following theorem establishes the asymptotic normality of the MELEs $(\hat\G_0,\hat\G_1)$.
\begin{theorem}
\label{thm1}
Assume that $\nu_i^*\in(0,1)$ for $i  = 0,1$ and Conditions C1--C4 are satisfied.
As the total sample size $n\to\infty$,
\begin{equation*}
n^{1/2}\ba{c}\hat\G_0 - \G_0^* \\ \hat\G_1 - \G_1^* \ea
\to
N\left({\bf 0}, \bSigma \right)
\end{equation*}
in distribution with the asymptotic variance-covariance matrix 
\begin{eqnarray}
\label{final_bSigma}
\bSigma = \frac{1}{\Delta^*}\J \Bigg[ E_0\left\{\frac{\bu(X)\bu(X)^\top}{h(X)}\right\} + \frac{1}{(\rho^*)^2}\B \Bigg] \J^\top
+ {\rm\diag}\lb \frac{\nu_0^*(1 - \G_0^*)^2}{\Delta^*(1-\rho^*)},\frac{\nu_1^*(1 - \G_1^*)^2}{\Delta^*\rho^*}\rb,
\end{eqnarray}
where 
\begin{eqnarray*}
&&\bu(x) = \left(x,u_0(x), \omega(x)x,\omega(x)u_0(x)\right)^\top,~
\tilde\bu_0(x) = -\rho^*\left(x,u_0(x)\right)^\top,\\
&&\tilde\bu_1(x) = (1-\rho^*)\left(x,u_1(x)\right)^\top,~ \tilde\bu(x) = \left(\tilde\bu_0(x)^\top,\tilde\bu_1(x) ^\top\right)^\top,\\
&&\B = E_0\{h_1(X)\tilde{\bu}(X)\Q(X)^\top\}\A_{\btheta}^{-1} E_0\{h_1(X)\Q(X)\tilde{\bu}(X)^\top\}.
\end{eqnarray*}
\end{theorem}

Note that the condition $\nu^*_i \in (0,1)$ for $i = 0,1$  ensures that the binomial log-likelihood function $\ell_0(\bnu)$ has regular properties
and the quadratic approximation is applicable. 
Repeating all the steps in the proof of Theorem \ref{thm1}, 
we obtain a similar result for cases where $\nu^*_i=0$ for $i=0,1$ in the following theorem.

\begin{theorem}
\label{thm2}
Assume that Conditions C1--C4 are satisfied.
When there is no excess of zeros, i.e., $\nu_i^* = 0$ for $i = 0,1$,  the joint distribution of $\sqrt{n}(\hat\G_0 - \G_0^*)$ and $\sqrt{n}(\hat\G_1 - \G_1^*)$ asymptotically follows a bivariate normal distribution with mean zero and  variance  in \eqref{final_bSigma} with $\nu_i^*$ being replaced by 0. 
\end{theorem}

Since the proposed method utilizes more information to obtain the MELEs of Gini indices, 
we expect that the proposed MELEs are more efficient than fully nonparametric estimators.
With the alternative form of the Gini index in \eqref{gini2}, 
the fully nonparametric estimators of the two Gini indices for sample $i=0,1$ are
\begin{eqnarray}
\label{tildeGini}
\tilde\G_i =(2\hat\nu_i-1)+(1-\hat\nu_i){\tilde\psi_i}/{\tilde m_i}, ~~i=0,1,
\end{eqnarray}
where 
\begin{eqnarray*}
\tilde m_i=  {n_{i1}}^{-1}\sum_{j=1}^{n_{i1}}X_{ij}, ~~
\tilde \psi_i=  \int_0^{\infty} \{2x \tilde G_i(x)\}d \tilde G_i(x),
\end{eqnarray*}
and 
$\tilde{G}_i(x) = {n_{i1}}^{-1}\sum_{j=1}^{n_{i1} }I(X_{ij} \leq x)$ 
is  the empirical CDF of the positive observations in sample $i$. 
The following theorem compares the proposed estimators MELEs $\hat\G_i$ and the nonparametric estimators $\tilde{\G}_i$ in terms of their asymptotic variance-covariance matrices.  Note that $\bSigma$ is the asymptotic variance-covariance matrix of the MELEs given in Theorem \ref{thm1}. 

\begin{theorem}
\label{thm3}
Assume that Conditions C1--C4 are satisfied. 
\begin{itemize}
    \item [(a)] For the nonparametric estimators  $(\tilde{\G}_0,\tilde{\G}_1)$ and as $n\to\infty$, we have 
\begin{equation*}
  \sqrt{n}
\left(
  \begin{array}{c}
\tilde\G_0-\G_0^*  \\
\tilde\G_1-\G_1^*\\
    \end{array}
  \right)
 \to N({\bf 0}, \bSigma_{non})  
\end{equation*}
in distribution, where the variance-covariance matrix 
\begin{eqnarray*}
 \bSigma_{non} &=&
 \J\diag\lb\frac{E_0\{\tilde{\bu}_0(X)\tilde{\bu}_0(X)^\top\}}{\Delta^*(\rho^*)^2(1-\rho)},\frac{E_0\{\omega(X)\tilde{\bu}_1(X)\tilde{\bu}_1(X)^\top\}}{\Delta^*\rho^*(1-\rho)^2}\rb\J^\top \\
 &&+ \diag\lb \frac{\nu_0^*(1 - \G_0^*)^2}{\Delta^*(1-\rho^*)},\frac{\nu_1^*(1 - \G_1^*)^2}{\Delta^*\rho^*}\rb .
\end{eqnarray*}

\item [(b)]
The two asymptotic variance-covariance matrices $\bSigma_{non}$ and $\bSigma$ satisfy 
\begin{equation*}
    \bSigma_{non}-\bSigma = \frac{1}{\Delta^*(\rho^*)^2(1-\rho^*)}\J  E_0\{h_1(X)\D(X)\D(X)^\top\}  \J^\top \geq {\bf 0},
\end{equation*}
where $\D(x) = \left(\D_0(x)^\top,\D_1(x)^\top\right)^\top$ for $x>0$ and 
\begin{equation*}
   \D_i(x)=\tilde{\bu}_i(x)- \Delta^*(1-\rho^*) E_0\left\{h_1(X)\tilde{\bu}_i(X)\Q(X)^\top\right\} \A_{\bstheta}^{-1} \Q(x),~i =0,1.
\end{equation*}
\end{itemize}
\end{theorem}

Note that $\bSigma_{non}-\bSigma$ is a positive semidefinite matrix, which implies that the proposed MELEs for the Gini index are at least as efficient as the nonparametric estimators. Our simulation results reported in Section~\ref{simu} confirm this result. 
It is worth mentioning that the theorem is applicable whether or not there are excess zero values.

\subsection{Inference on functions of Gini indices}

Under the current setting of two samples, we may be interested in performing inference on the Gini index for only one of the samples or other functions of the two Gini indices, such as their difference.
The results of Theorems \ref{thm1} and \ref{thm2} can be used to develop the following theorem  for parameters which are a general function of the two Gini indices. 

\begin{theorem}
\label{thm4}
Assume the conditions of Theorem \ref{thm3} hold. Let $\phi(\cdot,\cdot)$ be a bivariate smooth function. 
As $n \to \infty$, we have 
$\sqrt{n}\{\phi(\hat\G_0,\hat\G_1) -\phi(\G_0^*,\G_1^*)\} \to N(0, \sigma^2_{\phi})$ in distribution with 
$$
\sigma^2_{\phi} =\left( \frac{\partial \phi(\G_0^*,\G_1^*)}{\partial \G_0}, \frac{\partial \phi(\G_0^*,\G_1^*)}{\partial \G_1}\right) 
 \bSigma 
  \left( \frac{\partial \phi(\G_0^*,\G_1^*)}{\partial \G_0}, \frac{\partial \phi(\G_0^*,\G_1^*)}{\partial \G_1}\right) ^\top.
$$
\end{theorem}

With the results in Theorems \ref{thm3} and \ref{thm4}, we can  easily show that $\sigma^2_{\phi}$
is no larger than the asymptotic variance of the fully nonparametric estimator $\phi(\tilde\G_0,\tilde\G_1)$. 
That is, utilizing the information from both samples via the DRM \eqref{drm}
improves the estimation of $\phi(\G_0,\G_1)$. 

The general form $\phi(\cdot,\cdot)$ covers many interesting functions of $\G_0$ and $\G_1$. 
 For example,  when $\phi(x_1,x_2) = \mbox{logit}(x_1) = \log\{x_1/(1-x_1)\}$, the parameter $\phi(\G_0,\G_1)$ represents the logit transformation of the Gini index $\G_0$; 
 when $\phi(x_1,x_2) = x_1 - x_2$, the parameter 
$\phi(\G_0,\G_1)$ refers to the difference of two Gini  indices. 

The variance $\sigma^2_{\phi}$ may depend on $\G_0^*$, $\G_1^*$, and $(\boldsymbol{\nu},\boldsymbol{\theta}, \boldsymbol{P})$. 
Replacing these unknown quantities by their MELEs leads to a consistent estimator  $\hat\sigma^2_{\phi}$ of $\sigma_\phi^2$. 
Together with the result in Theorem \ref{thm4},  we have
$$
\sqrt{n}\{\phi(\hat\G_0,\hat\G_1) -\phi(\G_0^*,\G_1^*)\}/ \hat\sigma_{\phi}
\to N(0,1)
$$
in distribution. Hence, $\sqrt{n}\{\phi(\hat\G_0,\hat\G_1) -\phi(\G_0^*,\G_1^*)\}/ \hat\sigma_{\phi}$ 
is asymptotically pivotal and can be used to construct CIs  and to conduct hypothesis tests on $\phi(\G_0,\G_1)$. 

For ease of presentation, we use $\hat\phi$ and $\phi$ to denote  $\phi(\hat\G_0,\hat\G_1)$ and $\phi(\G_0,\G_1)$.
Then the $100(1-\tau)\%$ Wald-type CI for $\phi$ is given by 
$$
 [\hat\phi - z_{\tau/2}\hat\sigma_{\phi}/\sqrt{n},\hat\phi + z_{\tau/2}\hat\sigma_{\phi}/\sqrt{n}],
$$
where $z_{\tau/2}$  is the $(1-\tau/2)$ quantile of the standard normal distribution.
When testing $H_0: \phi = 0$, we reject the null hypothesis if 
$
|\sqrt{n}\hat\phi/\hat\sigma_{\phi}| > z_{\tau/2}
$
at the significance level $\tau$.

\section{Simulation Studies}
\label{simu}

In this section, we compare the finite-sample performance of our semiparametric methods with existing methods of inferences on the Gini indices through simulation studies. 
We focus on three inferential problems:
\begin{itemize}
    \item[(1)] Point estimation for  $\G_0$, $\G_1$, and  $\G_0 - \G_1$;
    \item[(2)] Confidence intervals on $\G_0$, $\G_1$, and  $\G_0 - \G_1$;
    \item[(3)] Hypothesis testing  on $H_0:\G_0=\G_1$. 
\end{itemize}

We conduct the simulation studies under two distributional settings:
(i) $G_0$ and $G_1$ are the CDFs of  $\chi^2_3$ and $\chi^2_4$;
and (ii) $G_0$ and $G_1$ are the CDFs of $Exp(0.5)$ and $Exp(1)$. 
Here $\chi^2_k$ represents the chi-square distribution with $k$ degrees of freedom, and $Exp(k)$ refers to the exponential distribution with the rate parameter $k$.
The proposed inference procedures under the DRM are implemented with the correctly specified $\q(x)$, where $\q(x) = \log(x)$ in the $\chi^2$ setting and $\q(x) = x$ in the exponential setting.
For each scenario, we consider two combinations of sample sizes, $(n_0,n_1) = (100,100)$, $(300,300)$, and the results are based on 2,000 Monte Carlo simulation runs.

\subsection{Performance of point estimators}
\label{simu_point}
We start by exploring the performance of the point estimators.  
We consider the following three estimators:
\begin{itemize}
    \item [--] \textit{EMP:}  $\tilde{\G}_0$, $\tilde{\G}_1$, and $\tilde{\G}_0 - \tilde{\G}_1$, where $\tilde{\G}_i$ is the nonparametric estimator given in \eqref{tildeGini} for $i = 0,1$;
    \item [--] \textit{JEL:} $\bar{\G}_0$, $\bar{\G}_1$, and $\bar{\G}_0 - \bar{\G}_1$, which are the jackknife empirical likelihood (JEL) estimators defined in \cite{wang2016jackknife}, where
    $$
    \bar{\G}_i = (2\tilde{\mu}_i)^{-1}\binom{n}{2}^{-1} \sum_{1\leq j_1 < j_2 \leq n_i}|X_{ij_1} - X_{ij_2}|
    $$ 
    with $\tilde{\mu}_i=\sum_{j=1}^{n_i}X_{ij}/n_i$ for $i = 0,1$;
    \item [--] \textit{DRM:} $\hat\G_0$, $\hat\G_1$, and $\hat\G_0-\hat\G_1$, where $\hat{\G}_i$ is the MELE given in \eqref{hatGini} for $i = 0,1$. 
\end{itemize}

Three combinations of $\bnu$ are considered for the zero population proportions: $(0,0)$, $(0.3,0.3)$, $(0.7,0.7)$. 
We evaluate the performance of a point estimator in terms of the bias and the mean squared error (MSE). 
Tables~\ref{point_chi} and \ref{point_exp} present the simulated results for different settings. 

\begin{table}[!htt]
  \centering
  \footnotesize	
  \tabcolsep 2mm
  \caption{Bias ($\times 1000$) and MSE ($\times 1000$) for point estimators ($\chi^2$).}
  \label{point_chi}
    \begin{tabular}{ccccccccc}
    \hline
          &       &       & \multicolumn{2}{c}{$\G_0$} & \multicolumn{2}{c}{$\G_1$} & \multicolumn{2}{c}{$\G_0-\G_1$} \\
          \cline{4-9}
    $(n_0,n_1)$ & $\bnu$   &       & Bias  & MSE   & Bias  & MSE   & Bias  & MSE \\
    \hline
    (100,100) & (0,0) & EMP   & 5.37  & 0.74  & 7.80  & 0.62  & -2.43 & 1.29 \\
          &       & JEL   & -0.39 & 0.73  & 1.57  & 0.57  & -1.96 & 1.31 \\
          &       & DRM   & 2.06  & 0.37  & 4.18  & 0.40  & -2.13 & 0.31 \\[1mm]
          & (0.3,0.3) & EMP   & 6.17  & 1.19  & 6.19  & 1.21  & -0.02 & 2.35 \\
          &       & JEL   & 2.16  & 1.18  & 1.83  & 1.20  & 0.33  & 2.40 \\
          &       & DRM   & 2.60  & 0.95  & 3.04  & 1.06  & -0.44 & 1.71 \\[1mm]
          & (0.7,0.7) & EMP   & 6.56  & 0.91  & 6.01  & 1.01  & 0.54  & 1.86 \\
          &       & JEL   & 4.88  & 0.91  & 4.18  & 1.01  & 0.70  & 1.89 \\
          &       & DRM   & 2.70  & 0.79  & 2.97  & 0.91  & -0.28 & 1.56 \\[2mm]
    (300,300) & (0,0) & EMP   & 1.70  & 0.23  & 2.31  & 0.19  & -0.61 & 0.40 \\
          &       & JEL   & -0.22 & 0.23  & 0.23  & 0.19  & -0.45 & 0.40 \\
          &       & DRM   & 0.66  & 0.13  & 1.12  & 0.14  & -0.46 & 0.11 \\[1mm]
          & (0.3,0.3) & EMP   & 2.52  & 0.39  & 1.75  & 0.41  & 0.78  & 0.79 \\
          &       & JEL   & 1.18  & 0.39  & 0.29  & 0.41  & 0.90  & 0.80 \\
          &       & DRM   & 0.94  & 0.32  & 0.89  & 0.37  & 0.06  & 0.57 \\[1mm]
          & (0.7,0.7) & EMP   & 2.84  & 0.31  & 2.40  & 0.34  & 0.44  & 0.67 \\
          &       & JEL   & 2.28  & 0.31  & 1.78  & 0.34  & 0.49  & 0.67 \\
          &       & DRM   & 1.30  & 0.27  & 1.50  & 0.31  & -0.20 & 0.55 \\
          \hline
    \end{tabular}%
\end{table}%

\begin{table}[!htt]
  \centering
  \footnotesize	
  \tabcolsep 2mm
  \caption{Bias ($\times 1000$) and MSE ($\times 1000$) for point estimators ($Exp$).}
  \label{point_exp}
    \begin{tabular}{ccccccccc}
    \hline
          &       &       & \multicolumn{2}{c}{$\G_0$} & \multicolumn{2}{c}{$\G_1$} & \multicolumn{2}{c}{$\G_0-\G_1$} \\
          \cline{4-9}
    $(n_0,n_1)$ & $\bnu$   &       & Bias  & MSE   & Bias  & MSE   & Bias  & MSE \\
    \hline
    (100,100) & (0,0) & EMP   & 4.59  & 0.81  & 5.94  & 0.86  & -1.35 & 1.59 \\
          &       & JEL   & -0.42 & 0.81  & 0.95  & 0.84  & -1.37 & 1.63 \\
          &       & DRM   & 1.55  & 0.66  & 2.82  & 0.41  & -1.27 & 0.58 \\[1mm]
          & (0.3,0.3) & EMP   & 4.85  & 1.09  & 3.53  & 1.10  & 1.31  & 2.14 \\
          &       & JEL   & 1.36  & 1.09  & 0.03  & 1.11  & 1.33  & 2.19 \\
          &       & DRM   & 1.64  & 0.96  & 0.73  & 0.82  & 0.91  & 1.46 \\[1mm]
          & (0.7,0.7) & EMP   & 5.17  & 0.81  & 3.62  & 0.80  & 1.55  & 1.55 \\
          &       & JEL   & 3.71  & 0.81  & 2.14  & 0.80  & 1.57  & 1.58 \\
          &       & DRM   & 1.81  & 0.73  & 1.24  & 0.65  & 0.56  & 1.20 \\[2mm]
    (300,300) & (0,0) & EMP   & 1.78  & 0.28  & 1.97  & 0.27  & -0.18 & 0.57 \\
          &       & JEL   & 0.11  & 0.28  & 0.30  & 0.27  & -0.19 & 0.58 \\
          &       & DRM   & 0.76  & 0.22  & 0.96  & 0.13  & -0.20 & 0.21 \\[1mm]
          & (0.3,0.3) & EMP   & 1.80  & 0.37  & 1.48  & 0.36  & 0.32  & 0.73 \\
          &       & JEL   & 0.63  & 0.37  & 0.31  & 0.36  & 0.32  & 0.74 \\
          &       & DRM   & 0.74  & 0.34  & 0.73  & 0.27  & 0.01  & 0.51 \\[1mm]
          & (0.7,0.7) & EMP   & 1.81  & 0.26  & 1.98  & 0.26  & -0.16 & 0.53 \\
          &       & JEL   & 1.32  & 0.26  & 1.48  & 0.26  & -0.16 & 0.54 \\
          &       & DRM   & 0.80  & 0.24  & 0.90  & 0.22  & -0.10 & 0.43 \\
          \hline
    \end{tabular}%
\end{table}%

We observe  from Tables~\ref{point_chi} and \ref{point_exp} that the biases of the estimators of $\G_0$ and $\G_1$ are acceptable for all three methods under all scenarios. The EMP estimators $\tilde{\G}_0$ and $\tilde{\G}_1$ always give the largest biases. When the proportions of zero values are small, i.e., $\bnu = (0,0)$ or $(0.3,0.3)$, the biases of the JEL estimators $\bar{\G}_0$ and $\bar{\G}_1$ are the smallest.
The DRM estimators $\hat{\G}_0$ and $\hat{\G}_1$ have a clear advantage in terms of bias when $\bnu = (0.7,0.7)$. 
The performance of the EMP estimators $\tilde{\G}_0$ and $\tilde{\G}_1$ and the JEL estimators $\bar{\G_0}$ and $\bar{\G_1}$ is similar in terms of the MSE. The DRM estimators $\hat{\G}_0$ and $\hat{\G}_1$ give the smallest MSEs in all cases; this agrees with the result in Theorem \ref{thm3}. 
The MSEs of all the estimators decrease as $\bnu$ moves toward $(0,0)$ or the sample size increases. 

For the estimators of the difference $\G_0 - \G_1$, we find that the biases of all the estimators are relatively small in all cases. The biases of the DRM estimator $\hat{\G}_0 - \hat{\G}_1$ are usually the smallest. The MSEs of the EMP estimator for ${\G}_0 - {\G}_1$ and JEL estimator for ${\G}_0 - {\G}_1$ are very close, whereas the MSEs of the DRM estimator are significantly smaller than those of the other two estimators. For instance, the MSE of $\hat{\G}_0 - \hat{\G}_1$ is less than 25\% of the MSEs of $\tilde{\G}_0 - \tilde{\G}_1$ and $\bar{\G}_0 - \bar{\G}_1$ when the simulated samples come from $\chi^2$ distributions with $(n_0,n_1)= (100,100)$ and $\bnu=(0,0)$.

We conducted additional simulations with $\bnu = (0.1,0.3)$ and $(0.6,0.4)$; the results show similar patterns and are presented in the Supplementary Material.

\subsection{Performance of confidence intervals}
\label{simu_ci}

We examine and compare the performance of the following confidence intervals for the Gini indices in the simulation studies:
\begin{itemize}
    \item [--] \textit{NA-EMP:} Wald-type CIs based on the normal approximation \citep{qin2010empirical};
    \item [--] \textit{BT-EMP:} bootstrap-t CIs  \citep{qin2010empirical};
    \item [--] \textit{EL}: ELR-based CIs  \citep{qin2010empirical};
    \item [--] \textit{BT-EL}: bootstrap ELR-based CIs  \citep{qin2010empirical};
    \item [--] \textit{JEL:} jackknife ELR-based CIs  \citep{wang2016jackknife,wang2016jackcompare};
    \item [--] \textit{AJEL:} adjusted jackknife ELR-based CIs  \citep{wang2016jackknife,wang2016jackcompare};
    \item [--] \textit{NA-DRM:} Wald-type CIs based on the normal approximation under the DRM;
    \item [--] \textit{BT-DRM:} bootstrap-t CIs under the DRM.
\end{itemize}

The EL method, to our best knowledge, has not been used to construct CIs for the difference of two Gini indices in the existing literature. Hence, we consider all the methods, except for \textit{EL} and \textit{BT-EL}, in our comparisons of the CIs for the parameter $\G_0-\G_1$.
For those calibrated by the bootstrap method, we used 1,000 nonparametric bootstrap samples drawn from the original sample with replacement. 

Three combinations of $\bnu$ are considered for the zero population proportions: $(0,0)$, $(0.3,0.3)$, $(0.7,0.7)$. 
We evaluate the performance of a CI in terms of the coverage probability (CP) and the average length (AL). 
Tables~\ref{CI_chi} and \ref{CI_exp} contain the simulated results for the CIs of $\G_0$ and $\G_1$ under different settings. 
The simulated results for the CIs of $\G_0-\G_1$ are shown in Table~\ref{CI_diff}.

\begin{table}[!htt]
  \centering
  \footnotesize	
  \tabcolsep 2mm
  \caption{Coverage probability (CP\%) and average length (AL) of CIs ($\chi^2$).}
  \label{CI_chi}
    \begin{tabular}{cc cccccccc}
    \hline
          & \multicolumn{1}{r}{} & \multicolumn{4}{c}{(100,100)}       & \multicolumn{4}{c}{(300,300)} \\
    \cline{3-10}
          & \multicolumn{1}{c}{} & \multicolumn{2}{c}{$\G_0$} & \multicolumn{2}{c}{$\G_1$} & \multicolumn{2}{c}{$\G_0$} & \multicolumn{2}{c}{$\G_1$} \\
    $\bnu$    & \multicolumn{1}{c}{} & \multicolumn{1}{c}{CP} & \multicolumn{1}{c}{AL} & \multicolumn{1}{c}{CP} & \multicolumn{1}{c}{AL} & \multicolumn{1}{c}{CP} & \multicolumn{1}{c}{AL} & \multicolumn{1}{c}{CP} & \multicolumn{1}{c}{AL} \\
    \hline
    (0,0) & NA-EMP & 93.85 & 0.100 & 94.20 & 0.092 & 94.60 & 0.059 & 94.80 & 0.054 \\
          & BT-EMP & 94.10 & 0.103 & 94.75 & 0.094 & 94.85 & 0.059 & 95.05 & 0.054 \\
          & EL    & 93.85 & 0.100 & 94.20 & 0.091 & 94.55 & 0.059 & 94.80 & 0.054 \\
          & BT-EL & 94.45 & 0.103 & 95.10 & 0.095 & 94.90 & 0.059 & 94.95 & 0.054 \\
          & JEL   & 94.45 & 0.102 & 94.85 & 0.094 & 94.70 & 0.059 & 95.15 & 0.054 \\
          & AJEL  & 94.80 & 0.105 & 95.50 & 0.096 & 94.90 & 0.060 & 95.30 & 0.055 \\
          & NA-DRM & 95.25 & 0.074 & 94.65 & 0.078 & 94.70 & 0.043 & 94.70 & 0.045 \\
          & BT-DRM & 95.55 & 0.075 & 95.00 & 0.079 & 94.55 & 0.043 & 94.55 & 0.046 \\[1.5mm]
    (0.3,0.3) & NA-EMP & 93.80 & 0.132 & 93.65 & 0.134 & 94.60 & 0.077 & 94.05 & 0.079 \\
          & BT-EMP & 95.30 & 0.135 & 94.55 & 0.137 & 95.20 & 0.077 & 94.40 & 0.079 \\
          & EL    & 93.75 & 0.131 & 93.65 & 0.134 & 94.60 & 0.077 & 94.00 & 0.078 \\
          & BT-EL & 94.50 & 0.136 & 94.85 & 0.139 & 94.65 & 0.078 & 94.55 & 0.079 \\
          & JEL   & 94.45 & 0.137 & 93.80 & 0.141 & 94.50 & 0.078 & 94.55 & 0.080 \\
          & AJEL  & 95.35 & 0.141 & 94.20 & 0.144 & 94.80 & 0.079 & 94.80 & 0.081 \\
          & NA-DRM & 95.10 & 0.120 & 94.35 & 0.130 & 95.45 & 0.070 & 94.90 & 0.076 \\
          & BT-DRM & 95.75 & 0.121 & 94.65 & 0.130 & 95.25 & 0.070 & 94.65 & 0.075 \\[1.5mm]
    (0.7,0.7) & NA-EMP & 92.20 & 0.113 & 92.95 & 0.119 & 94.90 & 0.067 & 93.90 & 0.070 \\
          & BT-EMP & 96.75 & 0.122 & 96.55 & 0.128 & 96.30 & 0.068 & 95.40 & 0.072 \\
          & EL    & 92.35 & 0.111 & 92.90 & 0.117 & 95.15 & 0.067 & 93.75 & 0.070 \\
          & BT-EL & 94.70 & 0.120 & 95.30 & 0.127 & 95.75 & 0.069 & 94.55 & 0.072 \\
          & JEL   & 90.75 & 0.123 & 90.80 & 0.129 & 94.00 & 0.069 & 93.00 & 0.072 \\
          & AJEL  & 91.35 & 0.127 & 91.55 & 0.133 & 94.25 & 0.070 & 93.10 & 0.073 \\
          & NA-DRM & 94.50 & 0.111 & 94.85 & 0.121 & 95.10 & 0.065 & 95.20 & 0.071 \\
          & BT-DRM & 95.40 & 0.113 & 95.90 & 0.123 & 95.45 & 0.064 & 95.60 & 0.070 \\
    \hline
    \end{tabular}%
\end{table}%

\begin{table}[!htt]
  \centering
  \footnotesize	
  \tabcolsep 2mm
  \caption{Coverage probability (CP\%) and average length (AL) of CIs ($Exp$).}
  \label{CI_exp}
    \begin{tabular}{cc cccccccc}
    \hline
          & \multicolumn{1}{r}{} & \multicolumn{4}{c}{(100,100)}       & \multicolumn{4}{c}{(300,300)} \\
    \cline{3-10}
          & \multicolumn{1}{c}{} & \multicolumn{2}{c}{$\G_0$} & \multicolumn{2}{c}{$\G_1$} & \multicolumn{2}{c}{$\G_0$} & \multicolumn{2}{c}{$\G_1$} \\
    $\bnu$    & \multicolumn{1}{c}{} & \multicolumn{1}{c}{CP} & \multicolumn{1}{c}{AL} & \multicolumn{1}{c}{CP} & \multicolumn{1}{c}{AL} & \multicolumn{1}{c}{CP} & \multicolumn{1}{c}{AL} & \multicolumn{1}{c}{CP} & \multicolumn{1}{c}{AL} \\
    \hline
    (0,0) & NA-EMP & 93.85 & 0.110 & 93.50 & 0.111 & 94.65 & 0.065 & 94.45 & 0.065 \\
          & BT-EMP & 94.35 & 0.115 & 94.05 & 0.115 & 94.75 & 0.065 & 94.75 & 0.065 \\
          & EL    & 93.90 & 0.110 & 93.50 & 0.110 & 94.65 & 0.065 & 94.55 & 0.065 \\
          & BT-EL & 94.50 & 0.113 & 94.00 & 0.113 & 94.80 & 0.065 & 94.60 & 0.065 \\
          & JEL   & 94.35 & 0.113 & 93.90 & 0.113 & 94.90 & 0.065 & 94.55 & 0.065 \\
          & AJEL  & 94.95 & 0.115 & 94.35 & 0.116 & 95.10 & 0.066 & 94.75 & 0.066 \\
          & NA-DRM & 94.80 & 0.100 & 94.05 & 0.079 & 93.95 & 0.059 & 95.20 & 0.045 \\
          & BT-DRM & 94.45 & 0.104 & 94.75 & 0.079 & 93.65 & 0.060 & 94.95 & 0.045 \\[1.5mm]
    (0.3,0.3) & NA-EMP & 93.55 & 0.127 & 94.25 & 0.128 & 94.55 & 0.075 & 93.45 & 0.075 \\
          & BT-EMP & 95.35 & 0.132 & 94.90 & 0.132 & 94.70 & 0.075 & 93.80 & 0.075 \\
          & EL    & 93.60 & 0.126 & 94.10 & 0.127 & 94.70 & 0.075 & 93.30 & 0.075 \\
          & BT-EL & 94.75 & 0.131 & 94.85 & 0.132 & 95.05 & 0.076 & 93.70 & 0.076 \\
          & JEL   & 93.80 & 0.132 & 94.55 & 0.133 & 94.65 & 0.076 & 93.65 & 0.076 \\
          & AJEL  & 94.50 & 0.136 & 95.15 & 0.136 & 95.00 & 0.076 & 93.80 & 0.076 \\
          & NA-DRM & 95.55 & 0.124 & 94.95 & 0.112 & 95.15 & 0.073 & 94.60 & 0.065 \\
          & BT-DRM & 95.45 & 0.125 & 95.30 & 0.112 & 94.60 & 0.072 & 94.60 & 0.064 \\[1.5mm]
    (0.7,0.7) & NA-EMP & 91.40 & 0.104 & 92.15 & 0.105 & 94.60 & 0.062 & 94.55 & 0.062 \\
          & BT-EMP & 96.30 & 0.114 & 95.70 & 0.115 & 95.85 & 0.064 & 95.50 & 0.064 \\
          & EL    & 92.05 & 0.102 & 92.35 & 0.102 & 94.70 & 0.062 & 94.55 & 0.062 \\
          & BT-EL & 95.00 & 0.110 & 94.20 & 0.111 & 95.40 & 0.064 & 95.25 & 0.064 \\
          & JEL   & 90.40 & 0.114 & 91.00 & 0.114 & 94.30 & 0.064 & 93.95 & 0.064 \\
          & AJEL  & 90.85 & 0.117 & 91.60 & 0.118 & 94.40 & 0.064 & 94.05 & 0.064 \\
          & NA-DRM & 94.65 & 0.109 & 93.90 & 0.101 & 96.10 & 0.064 & 95.40 & 0.059 \\
          & BT-DRM & 95.80 & 0.109 & 95.55 & 0.102 & 95.65 & 0.062 & 95.95 & 0.058 \\
    \hline
    \end{tabular}%
\end{table}%

When the sample sizes are $(100,100)$, we can see from Tables~\ref{CI_chi} and \ref{CI_exp} that the NA-EMP and EL CIs for $\G_0$ and $\G_1$ tend to be narrow and have lower CPs, especially when the proportions of zero values are large, i.e., $\bnu = (0.7,0,7)$. With the help of bootstrap calibration, the BT-EMP and BT-EL CIs achieve better performance in terms of CP. However, when $\bnu=(0.7,0.7)$, the BT-EMP CIs have slight overcoverage with inflated ALs. 
The AJEL CIs always have the longest ALs, and the JEL CIs are only slightly shorter. 
Moreover, when $\bnu = (0,0)$ and $(0.3,0.3)$, the CPs of the JEL and AJEL CIs are close to the nominal level of 95\%. 
The JEL and AJEL CIs suffer from undercoverage when $\bnu = (0.7,0.7)$. 
The NA-DRM CIs have the shortest ALs, and their CPs are very close to the 95\% nominal level in all cases. 
This is  strong evidence that using DRMs improves the performance of the CIs. 
The bootstrap calibration does little to improve the CIs: the performances of the NA-DRM and BT-DRM CIs are similar. 

When the sample sizes increase to $(300,300)$, the performance of all the CIs becomes satisfactory in terms of CP. The NA-DRM and BT-DRM CIs always have the shortest ALs, and there is little variation among the ALs of the other CIs. 

Since the Gini index ranges from 0 to 1, a logit transformation may improve the performance of the CIs for $\G_0$ and $\G_1$ under the DRM. However, the results (reported in the Supplementary Material) show that the transformation does not provide any significant improvement.

\begin{table}[!htt]
  \centering
  \footnotesize	
  \tabcolsep 2mm
  \caption{Coverage probability (CP\%) and average length (AL) of CIs for $\G_0 - \G_1$.}
  \label{CI_diff}
    \begin{tabular}{cc cccccccc}
    \hline
          & \multicolumn{1}{r}{} & \multicolumn{4}{c}{(100,100)}     & \multicolumn{4}{c}{(300,300)} \\
          \cline{3-10}
          & \multicolumn{1}{c}{} & \multicolumn{2}{c}{$\chi^2$} & \multicolumn{2}{c}{$Exp$} & \multicolumn{2}{c}{$\chi^2$} & \multicolumn{2}{c}{$Exp$} \\
          & \multicolumn{1}{c}{} & \multicolumn{1}{c}{CP} & \multicolumn{1}{c}{AL} & \multicolumn{1}{c}{CP} & \multicolumn{1}{c}{AL} & \multicolumn{1}{c}{CP} & \multicolumn{1}{c}{AL} & \multicolumn{1}{c}{CP} & \multicolumn{1}{c}{AL} \\
    \hline
    (0,0) & NA-EMP & 92.69 & 0.137 & 94.90 & 0.157 & 94.65 & 0.079 & 95.15 & 0.092 \\
          & BT-EMP & 92.89 & 0.139 & 94.80 & 0.161 & 94.35 & 0.080 & 95.15 & 0.092 \\
          & JEL   & 93.84 & 0.142 & 95.80 & 0.164 & 95.00 & 0.081 & 95.55 & 0.093 \\
          & AJEL  & 94.54 & 0.144 & 96.05 & 0.166 & 95.20 & 0.081 & 95.75 & 0.094 \\
          & NA-DRM & 94.44 & 0.070 & 94.95 & 0.092 & 94.65 & 0.041 & 94.70 & 0.055 \\
          & BT-DRM & 92.94 & 0.069 & 94.95 & 0.092 & 93.80 & 0.041 & 94.65 & 0.054 \\[1.5mm]
    (0.3,0.3) & NA-EMP & 94.19 & 0.188 & 93.05 & 0.181 & 95.00 & 0.110 & 94.65 & 0.106 \\
          & BT-EMP & 94.54 & 0.191 & 93.45 & 0.184 & 95.00 & 0.110 & 94.65 & 0.106 \\
          & JEL   & 95.45 & 0.202 & 94.75 & 0.195 & 95.55 & 0.113 & 95.35 & 0.108 \\
          & AJEL  & 95.70 & 0.205 & 95.10 & 0.198 & 95.55 & 0.113 & 95.55 & 0.109 \\
          & NA-DRM & 94.24 & 0.165 & 94.30 & 0.149 & 95.00 & 0.096 & 95.05 & 0.087 \\
          & BT-DRM & 93.34 & 0.161 & 93.45 & 0.146 & 94.35 & 0.094 & 94.60 & 0.086 \\[1.5mm]
 
    (0.7,0.7) & NA-EMP & 93.20 & 0.164 & 94.29 & 0.148 & 94.05 & 0.097 & 94.90 & 0.088 \\
          & BT-EMP & 93.65 & 0.170 & 94.14 & 0.153 & 94.05 & 0.098 & 94.90 & 0.089 \\
          & JEL   & 96.65 & 0.188 & 97.65 & 0.175 & 95.20 & 0.101 & 95.95 & 0.092 \\
          & AJEL  & 96.90 & 0.192 & 98.00 & 0.179 & 95.50 & 0.102 & 96.05 & 0.093 \\
          & NA-DRM & 95.55 & 0.162 & 96.19 & 0.138 & 95.60 & 0.094 & 95.70 & 0.080 \\
          & BT-DRM & 93.40 & 0.153 & 94.49 & 0.133 & 94.60 & 0.090 & 95.10 & 0.078 \\
    \hline
    \end{tabular}%
\end{table}%

We now discuss the simulation results for the CIs of the difference $\G_0-\G_1$ presented in Table \ref{CI_diff}. 
We observe that the NA-EMP and BT-EMP CIs have similar performance; their performance is acceptable except when the simulated samples are from $\chi^2$ distributions with $(n_0,n_1) = (100,100)$ and $\bnu = (0,0)$. In this case, the CPs of the NA-EMP and BT-EMP CIs are below the 95\% nominal level. 
The JEL and AJEL CIs always have the longest ALs. They experience overcoverage in some cases, especially when the proportions of zero values are high. 
The BT-DRM CIs have the shortest ALs, which leads to undercoverage in some cases. 
The performance of the NA-DRM CIs is consistently satisfactory in terms of CP and AL.

We also conduct additional simulations with $\bnu = (0.1,0.3)$ and $(0.6,0.4)$; the results display similar patterns and are presented in the Supplementary Material. 

\subsection{Performance of tests on the equality of two Gini indices}

In this section, 
we examine the performance of our proposed semiparametric test for testing the equality of the two Gini indices, i.e., $H_0:\G_0=\G_1$, with comparisons to other existing methods. 
We consider the following tests:
\begin{itemize}
    \item [--] \textit{NA-EMP:} Wald-type test based on the normal approximation of $ \tilde{\G}_0 - \tilde{\G}_1$ \citep{qin2010empirical};
    \item [--] \textit{NL-EMP:} Wald-type test based on the normal approximation of $ \mbox{logit}(\tilde{\G}_0) - \mbox{logit}(\tilde{\G}_1)$;
    \item [--] \textit{JEL:} jackknife ELR test \citep{wang2016jackcompare};
    \item [--] \textit{AJEL:} adjusted jackknife ELR test \citep{,wang2016jackcompare};
    \item [--] \textit{NA-DRM:} Wald-type test based on the normal approximation of $\hat{\G}_0 - \hat{\G}_1$ under the DRM;
    \item [--] \textit{NL-DRM:} Wald-type test based on the normal approximation of $\mbox{logit}(\hat{\G}_0) -\mbox{logit}(\hat{\G}_1)$ under the DRM.
\end{itemize}

Several combinations of $\bnu$ are chosen to satisfy the null hypothesis $H_0$ or the alternative  hypothesis $H_a$. The details are presented in Table \ref{choice_bnu}. 
Tables \ref{simu_typeone} and \ref{simu_testpower} give the simulated type I error rate and simulated  power of each test at the 5\% significance level. 

\begin{table}[!htt]
  \centering
  \footnotesize	
  \caption{Choices of $\bnu$ in simulations of testing the equality of the two Gini indices.}
  \label{choice_bnu}
    \begin{tabular}{ccccccc}
    \hline
    &\multicolumn{6}{c}{Null hypothesis $H_0$}\\
     \cline{2-7}
    & \multicolumn{3}{c}{$\chi^2$} & \multicolumn{3}{c}{$Exp$} \\
    $\bnu$    & (0,0.079) & (0.3,0.355) & (0.7,0.724) & (0,0) & (0.3,0.3) & (0.7,0.7) \\
    \hline 
    &\multicolumn{6}{c}{Alternative hypothesis $H_a$}\\
    \cline{2-7}
    & \multicolumn{3}{c}{$\chi^2$} & \multicolumn{3}{c}{$Exp$} \\
  $\bnu$    & (0,0) & (0.1,0.3) & (0.4,0.65) & (0.1,0.3) & (0.3,0.45) & (0.5,0.4) \\
  $\G_0 - \G_1$  & 0.049 & -0.081 & -0.127 & -0.100 & -0.075 & 0.050 \\
  $\mbox{logit}(\G_0) -\mbox{logit}(\G_1)$ & 0.206 & -0.323 & -0.633 & -0.418 & -0.350 & 0.251 \\
          \hline
    \end{tabular}%
\end{table}%

\begin{table}[!htt]
  \centering
  \footnotesize	
  \tabcolsep 2mm
  \caption{Type I error rate (\%) for testing  $H_0:\G_0=\G_1$ at the 5\% significance level.}
  \label{simu_typeone}
    \begin{tabular}{cccccccc}
    \hline
          &       & \multicolumn{3}{c}{$\chi^2$} & \multicolumn{3}{c}{$Exp$} \\
          \cline{3-8}
          & $\bnu$    & \multicolumn{1}{c}{(0,0.079)} & \multicolumn{1}{c}{(0.3,0.355)} & \multicolumn{1}{c}{(0.7,0.724)} & \multicolumn{1}{c}{(0,0)} & \multicolumn{1}{c}{(0.3,0.3)} & \multicolumn{1}{c}{(0.7,0.7)} \\
          \hline
    (100,100) & NA-EMP & 5.50  & 5.00  & 7.10  & 5.55  & 4.85  & 6.90 \\
          & NL-EMP & 5.45  & 4.85  & 6.40  & 5.50  & 4.70  & 6.30 \\
          & JEL   & 4.55  & 4.05  & 3.70  & 4.85  & 3.10  & 2.70 \\
          & AJEL  & 4.15  & 3.85  & 3.45  & 4.45  & 2.70  & 2.40 \\
          & NA-DRM & 4.90  & 5.15  & 5.15  & 5.05  & 4.70  & 5.20 \\
          & NL-DRM & 4.85  & 4.80  & 4.75  & 4.95  & 4.65  & 4.95 \\ [1.5mm]
    (300,300) & NA-EMP & 5.10  & 5.70  & 5.70  & 6.15  & 5.35  & 5.55 \\
          & NL-EMP & 5.10  & 5.70  & 5.70  & 6.15  & 5.35  & 5.55 \\
          & JEL   & 4.90  & 5.35  & 5.10  & 5.70  & 4.85  & 4.20 \\
          & AJEL  & 4.80  & 5.20  & 4.80  & 5.55  & 4.80  & 4.00 \\
          & NA-DRM & 5.05  & 4.90  & 4.90  & 5.25  & 5.30  & 5.15 \\
          & NL-DRM & 5.05  & 4.85  & 4.95  & 5.25  & 5.25  & 5.05 \\
          \hline
    \end{tabular}%
\end{table}%

From Table~\ref{simu_typeone}, we observe that the type I error rates for NA-DRM are stable and close to the 5\% significance level in all cases. The type I error rates for NL-DRM are similar to those for NA-DRM when the sample sizes are $(300,300)$ and smaller when $(n_0,n_1) = (100,100)$. 
This implies that the logit transformation of the Gini indices is unnecessary for the equality test.
The type I error rates for NA-EMP and NL-EMP show similar trends. 
When the sample sizes are $(100,100)$ and the proportions of zero values are high,  NA-EMP, NL-EMP, JEL, and AJEL have either inflated  or conservative  type I error rates. Large sample sizes seem to improve their performance.

\begin{table}[!htt]
  \centering
  \footnotesize	
  \tabcolsep 2mm
  \caption{Simulated testing power (\%) of rejecting  $H_0: \G_0=\G_1$ at the 5\% significance level.}
  \label{simu_testpower}
    \begin{tabular}{cccccccc}
    \hline
          &       & \multicolumn{3}{c}{$\chi^2$} & \multicolumn{3}{c}{$Exp$} \\
          \cline{3-8}
          & $\bnu$    & (0,0) & (0.1,0.3) & (0.4,0.65) & (0.1,0.3) & (0.3,0.45) & (0.5,0.4) \\

          \hline
    (100,100) & NA-EMP & 30.15 & 43.45 & 78.05 & 59.65 & 38.20 & 18.50 \\
          & NL-EMP & 30.00 & 42.95 & 76.95 & 59.10 & 37.40 & 18.15 \\
          & JEL   & 28.85 & 42.90 & 75.70 & 58.00 & 33.75 & 14.90 \\
          & AJEL  & 28.00 & 42.05 & 74.85 & 56.95 & 32.75 & 14.00 \\
          & NA-DRM & 82.60 & 58.35 & 83.20 & 80.75 & 50.05 & 23.20 \\
          & NL-DRM & 82.45 & 58.05 & 82.25 & 80.45 & 49.95 & 22.05 \\[1.5mm]
    (300,300) & NA-EMP & 67.30 & 85.80 & 99.75 & 97.00 & 79.85 & 45.90 \\
          & NL-EMP & 67.30 & 85.70 & 99.70 & 96.90 & 79.50 & 45.50 \\
          & JEL   & 66.90 & 86.10 & 99.65 & 97.10 & 79.05 & 44.70 \\
          & AJEL  & 66.50 & 85.85 & 99.65 & 96.85 & 78.60 & 44.05 \\
          & NA-DRM & 99.95 & 95.70 & 99.85 & 99.90 & 90.75 & 56.90 \\
          & NL-DRM & 99.95 & 95.70 & 99.80 & 99.90 & 90.80 & 55.75 \\
          \hline
    \end{tabular}%
\end{table}%

We observe from Table~\ref{simu_testpower} that NA-DRM always gives the largest testing powers. The performance of NL-DRM is comparable to NA-DRM.  
When the true difference of the Gini indices is large, the testing powers of NA-DRM and NL-DRM are significantly larger than those of the other methods. For example, when the simulated samples are from the $\chi^2$ distributions with $(n_0,n_1) = (100,100)$ and $\bnu= (0,0)$, the testing powers of NA-DRM and NL-DRM are more than twice the others. 


\section{Real-Data Applications}
\label{realdata}

In this section, we apply our proposed methods to analyze two real datasets. 
Each dataset can be viewed as consisting of two samples from two different populations, and we are interested in computing the point estimates as well as the construction of $95\%$ confidence intervals for the Gini indices and their difference. The populations for the first dataset contain a large proportions of zeros and the study variables for the second dataset are strictly positive. 

The first dataset \citep{zhou2008computer} is from a clinical drug utilization study of patients with uncomplicated hypertension, originally conducted by  \cite{murray2004failure}. 
It consists of the inpatient charges of 483 patients by gender. 
We label the charges of the 282 male patients as sample 0 and those of the 201 female patients as sample 1. 
In most cases, uncomplicated hypertension can be controlled if the patients follow guidelines and take antihypertensive drugs regularly. If they do not need inpatient treatment, the corresponding charges are zero. 
There are 253 zero values (89.7\%) im sample 0 and 171 (85.0\%) in sample 1.  

\begin{figure}
    \centering
    \includegraphics[scale = 0.65]{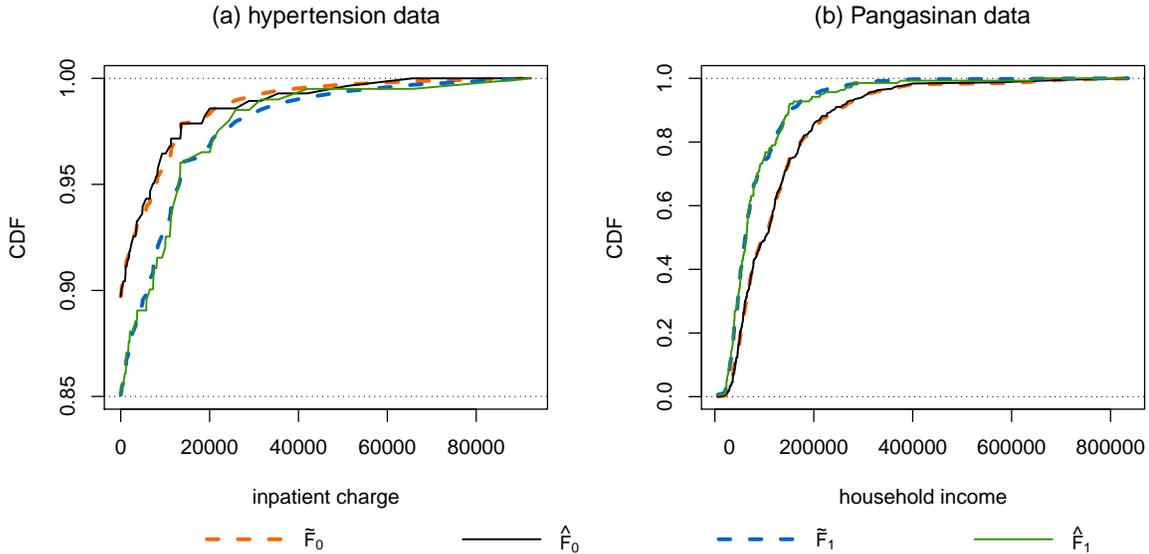}
    \caption{Fitted population distributions for real datasets. $\hat{F}_0$ and $\hat{F}_1$: fitted CDFs under the DRM; $\tilde{F}_0$ and $\tilde{F}_1$: empirical CDFs. }
    \label{cdf_plot}
\end{figure}

To analyze the dataset with our proposed method, we need to choose an appropriate $\q(x)$ in the DRM \eqref{drm}. The dataset is highly
skewed to the right because of the high proportions of zero values and extra skewness in the positive
inpatient charges. To balance model fit and model complexity, we choose $\q(x) = \log(x)$. The goodness-of-fit test of \cite{qin1997goodness} gives a p-value of 0.563, which indicates that this is a suitable choice. 
Figure \ref{cdf_plot}(a) shows the fitted population distribution functions $\hat F_0$ and $\hat F_1$ under the DRM with $\q(x) = \log(x)$ together with the empirical CDFs $\tilde{F}_0$ and $\tilde{F}_1$. 
Clearly, the fit is adequate.  

We apply the methods discussed in Sections~\ref{simu_point} and \ref{simu_ci} to this dataset.
Table~\ref{point_hypertension} presents the point estimates, and Table~\ref{ci_hypertension} shows the lower bound (LB), upper bound (UB),
 and length of the 95\% CIs. 
The estimates of $\G_0$, $\G_1$, and $\G_0-\G_1$ for all three methods are very close. In particular, the EMP and JEL estimates   are almost the same. 
The estimates of $\G_0$ and $\G_1$ are greater than 0.93, indicating the large inequality of the inpatient charge for patients with uncomplicated hypertension; 
the high proportion of zero values contributes to this. 
All the methods give similar 95\% CIs for $\G_0$. The 95\% CIs for $\G_1$ and $\G_0 - \G_1$ for NA-DRM and BT-DRM are the shortest. 
All the CIs for $\G_0 - \G_1$ contain 0, which suggests no significant difference between the inequality of the inpatient charge for female and male patients at the 95\% confidence level.

\begin{table}[!htt]
  \centering
  \footnotesize	
  \tabcolsep 2mm
  \caption{Point estimates of Gini indices and their difference (hypertension data).}
  \label{point_hypertension}
    \begin{tabular}{cccc}
    \hline
          & $\G_0$ (male)    & $\G_1$  (female)   & $\G_0-\G_1$ \\
          \hline
    EMP   & 0.959 & 0.933 & 0.026 \\
    JEL   & 0.959 & 0.933 & 0.026 \\
    DRM   & 0.956 & 0.934 & 0.022 \\
    \hline
    \end{tabular}%
\end{table}%

\begin{table}[!htt]
  \centering
  \footnotesize	
  \tabcolsep 2mm
  \caption{95\% CIs for the two Gini indices and their difference (hypertension data).}
  \label{ci_hypertension}
    \begin{tabular}{cccccccccc}
    \hline
          & \multicolumn{3}{c}{$\G_0$ (male)} & \multicolumn{3}{c}{$\G_1$ (female)} & \multicolumn{3}{c}{$\G_0- \G_1$} \\
          \cline{2-10}
          & LB    & UB    & Length & LB    & UB    & Length & LB    & UB    & Length \\
          \hline
    NA-EMP & 0.942 & 0.977 & 0.035 & 0.902 & 0.964 & 0.062 & -0.009 & 0.062 & 0.071 \\
    BT-EMP & 0.936 & 0.974 & 0.039 & 0.888 & 0.959 & 0.071 & -0.013 & 0.066 & 0.078 \\
    EL    & 0.941 & 0.975 & 0.034 & 0.903 & 0.961 & 0.058 &  --     &   --    &  --\\
    BT-EL & 0.938 & 0.976 & 0.038 & 0.897 & 0.966 & 0.068 &  --     &  --     & -- \\
    JEL   & 0.942 & 0.980 & 0.038 & 0.904 & 0.967 & 0.063 & -0.017 & 0.069 & 0.086 \\
    AJEL  & 0.942 & 0.981 & 0.039 & 0.904 & 0.967 & 0.064 & -0.018 & 0.069 & 0.087 \\
    NA-DRM & 0.938 & 0.974 & 0.036 & 0.906 & 0.961 & 0.056 & -0.009 & 0.054 & 0.063 \\
    BT-DRM & 0.934 & 0.972 & 0.038 & 0.901 & 0.957 & 0.056 & -0.007 & 0.048 & 0.055 \\
    \hline
    \end{tabular}%
\end{table}%

The second dataset comes from the 1997 Family and Income and Expenditure Survey conducted by the Philippine Statistics Authority; the metadata is available in the R package \texttt{ineq}.
The province of Pangasinan is located in the Ilocos Region of Luzon. 
The dataset contains household incomes from different areas of Pangasinan: urban (Sample 0) and rural (Sample 1).
Sample 0 has 245 observations and sample 1 has 138 observations. All the incomes are positive.

The skewness of the dataset suggests setting $\q(x) = \log(x)$ in the DRM \eqref{drm}. The goodness-of-fit test of \cite{qin1997goodness} gives a p-value 0.607. Hence, there is no strong evidence to reject the choice of  $\q(x) = \log(x)$. Figure~\ref{cdf_plot}(b) also shows that the DRM with $\q(x) = \log(x)$ fits the data well. 

We use all the methods of Sections~\ref{simu_point} and \ref{simu_ci} to analyze the dataset and summarize the results in Tables~\ref{point_Pangasinan} and \ref{ci_Pangasinan}.
The EMP and JEL methods give similar estimates of $\G_0$, $\G_1$, and $\G_0-\G_1$. 
The DRM estimate of $\G_0$ is comparable to the other estimates, while the DRM estimate of $\G_1$ is smaller than the others. 
Hence, the DRM estimate of $\G_0 - \G_1$ is larger. 
All the methods give similar results for the 95\% CIs for $\G_0$. The 95\% CIs for $\G_1$ and $\G_0 - \G_1$ by NA-DRM and BT-DRM are significantly shorter than the other CIs. 
This is strong evidence that our method helps to utilize information across the two samples and effectively improves inference when sample sizes are small or moderate. 
We do not reject the hypothesis that the income inequalities of urban and rural households are the same, since all the 95\% CIs for $\G_0 - \G_1$ contain 0.

\begin{table}[!htt]
  \centering
  \footnotesize	
  \tabcolsep 2mm
  \caption{Point estimates of Gini indices and their difference (Pangasinan data).}
  \label{point_Pangasinan}
    \begin{tabular}{cccc}
    \hline
      & $\G_0$ (urban)    & $\G_1$ (rural)   & $\G_0-\G_1$ \\
      \hline
    EMP   & 0.393 & 0.394 & -0.001 \\
    JEL   & 0.391 & 0.389 & 0.002 \\
    DRM   & 0.399 & 0.371 & 0.028 \\
    \hline
    \end{tabular}%
\end{table}%

\begin{table}[!htt]
  \centering
  \footnotesize	
  \tabcolsep 2mm
  \caption{95\% CIs for the two Gini indices and their difference (Pangasinan data).}
  \label{ci_Pangasinan}
    \begin{tabular}{cccccccccc}
     \hline
          & \multicolumn{3}{c}{$\G_0$ (urban) } & \multicolumn{3}{c}{$\G_1$ (rural)} & \multicolumn{3}{c}{$\G_0- \G_1$} \\
          \cline{2-10}
          & LB    & UB    & Length & LB    & UB    & Length & LB    & UB    & Length \\
          \hline
    NA-EMP & 0.354 & 0.433 & 0.079 & 0.332 & 0.455 & 0.123 & -0.074 & 0.073 & 0.146 \\
    BT-EMP & 0.356 & 0.441 & 0.085 & 0.338 & 0.481 & 0.143 & -0.085 & 0.068 & 0.153 \\
    EL    & 0.354 & 0.433 & 0.079 & 0.333 & 0.456 & 0.123 &   --    &    --   & -- \\
    BT-EL & 0.353 & 0.434 & 0.080 & 0.335 & 0.455 & 0.120 &   --    &  --     &  -- \\
    JEL   & 0.356 & 0.436 & 0.081 & 0.339 & 0.466 & 0.127 & -0.083 & 0.070 & 0.153 \\
    AJEL  & 0.355 & 0.437 & 0.081 & 0.338 & 0.467 & 0.129 & -0.084 & 0.071 & 0.154 \\
    NA-DRM & 0.361 & 0.436 & 0.075 & 0.343 & 0.399 & 0.055 & -0.003 & 0.059 & 0.062 \\
    BT-DRM & 0.359 & 0.443 & 0.084 & 0.343 & 0.403 & 0.060 & -0.006 & 0.057 & 0.063 \\
    \hline
    \end{tabular}%
\end{table}%

\section{Concluding Remarks}
\label{conclude}

We have proposed new semiparametric inference procedures for the Gini indices of two semicontinuous populations.
Under the mixture model \eqref{mixture} and the DRM \eqref{drm},  
we proposed the MELEs of the Gini indices $\G_0$ and $\G_1$ and established
the asymptotic normality of the MELEs.
Our methods are applicable whether or not there are excess zero values. 
We showed numerically and theoretically that our MELEs are more efficient than fully nonparametric estimators.
We also explored the asymptotic properties of a general function of two Gini indices, and used the difference of the two Gini indices  as an illustrating example.
We used the asymptotic results to construct CIs and perform hypothesis tests for $\G_0$, $\G_1$, and $\G_0 - \G_1$.
Simulation studies demonstrated that the CIs under the DRM have superior performance in terms of coverage accuracy and average length. 
Moreover, our method has a higher testing power than existing methods.

ELR-based CIs are range-preserving. It would be interesting to construct ELR-based CIs for Gini indices under the DRM, but 
the theoretical development may be technically challenging.
\cite{qin2010empirical} investigated inference on the Gini index of a population under stratified random sampling. We could use the DRM to link the distributions of the subpopulations in all strata  and then develop an inference procedure for the Gini index of the whole population.
Correlated data is also of interest: \cite{peng2011empirical} and \cite{wang2016jackcompare} compared the Gini indices of paired data.
The DRM would be useful for modeling the marginal distributions  of paired data and improving the  efficiency of the estimation of the Gini indices. 
We leave these topics to future research.


\newpage

{\centering {\large {\bf Supplementary material for \\
``Semiparametric inference on Gini indices of two semicontinuous populations under density ratio models''}}\par }

 \setcounter{equation}{0}
\setcounter{section}{0}
\renewcommand{\theequation} {S.\arabic{equation}}
\setcounter{theorem}{0}
\setcounter{table}{0}

\bigskip

\noindent
This document of supplementary material provides further details for the paper entitled ``Semiparametric inference on Gini indices of two semicontinuous populations under density ratio models". 
It contains the proofs of Theorems 1--4 in the main paper and additional simulation results. 
The basic setting and some useful lemmas are presented in Section 1. 
The proofs and technical details for Theorems 1--4 are given in Sections 2--5. 
Section 6 contains some additional simulation results.

\setcounter{section}{0}
\section{The basic setting and useful lemmas}

Recall that 
\begin{equation}
 \label{App_mixture}
  (X_{i1},\cdots, X_{in_i}) \sim F_i(x)= \nu_iI(x\geq0)+(1-\nu_i)I(x>0)G_i(x), ~~\mbox{for~}~i=0,1,
\end{equation}
where $\nu_i$ is the proportion of zeros in sample $i$, $n_i$ is the sample size for sample $i$, $I(\cdot)$ is an indicator function, 
and  $G_i(\cdot)$ is the cumulative distribution function (CDF) of the positive observations in sample $i$.
We link $G_0(x)$ and $G_1(x)$ via a density ratio model (DRM):
\begin{equation}
 \label{App_drm}
dG_1(x) = \exp\{\alpha +\boldsymbol{\beta}^\top\q(x)\}dG_0(x) = \exp\{\btheta^\top\Q(x)\}dG_0(x)
\end{equation}
for the unknown parameter $\btheta = (\alpha,\bbeta^\top)^\top$ and $\Q(x) = (1,\q(x)^\top)^\top$ with $\q(x)$ a $d$-dimensional, prespecified, 
nontrivial basis function.

Recall that $n_{i0}$ and $n_{i1}$ are the (random) numbers of zero observations and positive observations, respectively, in sample $i=0,1$.  Clearly, $n_i = n_{i0}+n_{i1}$, for $i=0,1$.
Without loss of generality, we assume that the first $n_{i1}$ observations in group $i$, $X_{i1},\cdots, X_{in_{i1}}$, are positive,
and the remaining $n_{i0}$ observations are 0.
Let $n$ be the total (fixed) sample size, i.e., $n=n_0+n_1$.

We further let $\bnu = (\nu_0,\nu_1)^\top$.
The maximum empirical likelihood estimators (MELEs) of $\bnu$ and $\btheta$ respectively maximize
$
\ell_0\left(\bnu\right)
$
and
$
\ell_1(\btheta)
$, where
$$
\ell_0\left(\bnu\right) = \sum_{i=0}^1\log\left\{ \nu_i^{n_{i0}}\left( 1-\nu_i\right)^{n_{i1}} \right\}
$$
and
$$
\ell_1(\btheta)=
- \Sumij
  \log\left\{ 1 + \hat\rho[\exp\{\btheta^{\top} \Q(X_{ij})\}-1] \right\}
+ \sum_{j=1}^{n_{11}} \{\btheta^{\top}  \Q(X_{1j}) \}
$$
with $\hat\rho = {n_{11}}/{(n_{01}+n_{11})}$ being a random variable.
That is,
\begin{equation}
 \label{App_mele1}
\hat\bnu=\arg\max_{\bsnu} \ell_0\left(\bnu\right)
~~
\mbox{ and }
~~
\hat\btheta=\arg\max_{\bstheta}\ell_1(\btheta).
\end{equation}
Once $\hat {\boldsymbol{\theta}} $ is obtained,
we have 
$$\hat p_{ij} = (n_{01}+n_{11})^{-1} \lb 1 + \hat{\rho}[\exp\{ \hat\btheta^{\top} \Q(X_{ij})\}-1]\rb ^{-1}.$$
Note that $\sum_{i=0}^1\sum_{j=1}^{n_{i1}}\hat p_{ij} = 1$, 
which ensures that the MELE of $G_0(x)$ is a CDF.
The MELEs of $G_0(x)$ and $G_1(x)$ for $x >0$ are 
\begin{equation}
 \label{App_hatCDF}
    \hat G_0(x) = \Sumij \hat p_{ij} I(X_{ij}\leq x)~~
    \mbox{and}
    ~~
    \hat G_1(x) = \Sumij \hat p_{ij}\exp\{ \hat\btheta^{\top} \Q(X_{ij})\} I(X_{ij}\leq x).
\end{equation}

For convenience of presentation, we recall and introduce some notation.
 We use $\bnu^{*}$ and $\btheta^{*}$ to denote the true values of $\bnu$ and $\btheta$. 
Let $\etab =(\bnu^\top,\rho,\btheta^\top)^\top $, $w_i = n_i/n$ for $i = 0,1$, and
\begin{eqnarray*}
&\Delta^*=\sum_{i=0}^1 w_i(1-\nu_i^*), 
~\rho^*=\frac{w_1(1-\nu_1^*)}{\Delta^*},~
\omega(x;\btheta)=\exp\{\btheta^{\top} \Q(x)\},~\omega(x) = \omega(x;\btheta^*),\\
&h(x)=1+\rho^* \{\omega(x)-1\},~h_1(x)=\rho^* \omega(x)/h(x),~h_0(x) =(1-\rho^*)/h(x).
\end{eqnarray*}

Note that $\omega(\cdot)$, $h(\cdot)$, $h_0(\cdot)$, and $h_1(\cdot)$ depend on $\btheta^*$ and/or $\rho^*$ and $h_0(x) + h_1(x) =1$. Henceforth, we use $\sum_{ij}$ and $\sum_{ls}$ to denote summation over the full range of data. 

\subsection{Alternative form of Gini index}
According to \cite{david1968miscellanea}, the Gini's mean difference for sample $i$ can be equivalently expressed by 
\begin{equation*}
D_i= E|X_{i1}-X_{i2}| 
= 2\int_{-\infty}^{\infty} \left\{2xF_i(x) - x\right\} dF_i(x). 
\end{equation*}
Under model  \eqref{App_mixture}, $F_i(x)= \nu_iI(x \geq 0)+(1-\nu_i)I(x>0)G_i(x)$.
Then 
$D_i$ can be further written as 
\begin{eqnarray*}
    D_i &=& 2(1-\nu_i)\int_0^{\infty} [2x \{\nu_i + (1-\nu_i)G_i(x)\} -x]dG_i(x)\\
    &=& 2(1-\nu_i)\int_0^{\infty} x \{(2 \nu_i - 1) + (1-\nu_i)2G_i(x)\} dG_i(x)\\
    &=& 2(2 \nu_i - 1)\int_0^{\infty} x(1-\nu_i) dG_i(x) + 2(1-\nu_i)^2\int_0^{\infty}\{2xG_i(x)\} dG_i(x)\\
    &=&2(2 \nu_i - 1)\mu_i + 2(1-\nu_i)^2\int_0^{\infty}\{2xG_i(x)\} dG_i(x).
\end{eqnarray*}
Let $m_i = \int_0^{\infty} x dG_i(x) $ and $\psi_i = \int_0^{\infty} \{2xG_i(x)\}dG_i(x)$. We then have 
$\mu_i = (1-\nu_i)m_i$ and 
\begin{equation*}
   \G_i = \frac{D_i}{2\mu_i}
= (2 \nu_i - 1) + (1-\nu_i) \frac{\psi_i}{m_i}.
\end{equation*}

With the definition of $m_i$ and $\psi_i$, and the MELEs of CDFs $G_i$'s in  \eqref{App_hatCDF}, the MELEs of $m_i$ and $\psi_i$
are as follows: 
\begin{eqnarray*}
\hat m_0 &=& \sum_{i=0}^1\sum_{j=1}^{n_{i1}} \hat{p}_{ij} X_{ij},~~~
\hat m_1= \sum_{i=0}^1\sum_{j=1}^{n_{i1}} \hat{p}_{ij}\omega(X_{ij};\hat\btheta) X_{ij},\\
\hat \psi_0 &=& \sum_{i=0}^1\sum_{j=1}^{n_{i1}} \hat{p}_{ij} \{2X_{ij}\hat G_0(X_{ij})\}
 = \sum_{i =0}^1\sum_{j =1}^{n_{i1}} \hat{p}_{ij} X_{ij}\lb 2 \sum_{l =0}^1 \sum_{s = 1}^{n_{l1}} \hat{p}_{ls}I(X_{ls} \leq X_{ij}) \rb,\\
\hat \psi_1&=& \sum_{i=0}^1\sum_{j=1}^{n_{i1}} \hat{p}_{ij}\omega(X_{ij};\hat\btheta) \{2X_{ij}\hat G_1(X_{ij})\} \\
&&= \sum_{i =0}^1\sum_{j =1}^{n_{i1}} \hat{p}_{ij}\omega(X_{ij};\htheta) X_{ij}\lb 2 \sum_{l =0}^1 \sum_{s = 1}^{n_{l1}} \hat{p}_{ls}\omega(X_{ls};\htheta)I(X_{ls} \leq X_{ij}) \rb.
\end{eqnarray*}

The MELEs of the two Gini indices are given by
\begin{eqnarray}
 \label{App_Gini_alternative}
\hat\G_i=(2 \hat{\nu}_i - 1) + (1-\hat{\nu}_i) \frac{\hat\psi_i}{\hat m_i} ~~~\text{for}~i=0,1.
\end{eqnarray}

\subsection{Some useful lemmas}

We present several useful lemmas in preparation for the proofs in Sections 2--5. 
The first lemma considers the expectation of summations. 
\begin{lemma}
 \label{App_lem1}
Suppose that $g(x)$ is an arbitrary vector-valued function.
Let $E_0(\cdot)$ represent the expectation operator with respect to $G_0$ and $X$ be a random variable from $G$.
Then 
$$
E \left\{\sum_{ij}g(X_{ij}) I(X_{ij}>0) \right\} = \sum_{i=0}^1 n_{i}{ E}\{g(X_{i1}) I(X_{i1}>0) \} = n\Delta^* E_0\{h(X)g(X)\} .
$$
\end{lemma}

\proof
Under the DRM ( \ref{App_drm}),
\begin{eqnarray*}
\sum_{i=0}^1 n_{i}{ E}\{g(X_{i1}) I(X_{i1}>0) \} 
= n_{0}(1-\nu_0^*) E_0\{g(X) \} + n_{1}(1-\nu_1^*) E_0\{\omega(X)g(X) \}.
\end{eqnarray*}
Since $w_i=n_i/n$ and using the definitions of $\Delta^*$ and $\rho^*$, we further have 
\begin{eqnarray*}
E\left\{\sum_{ij} g(X_{ij}) I(X_{ij}>0) \right\}
&=&n w_0(1-\nu_0^*) E_0\{g(X) \} + nw_1(1-\nu_1^*) E_0\{\omega(X)g(X) \}\\
&=& n\Delta^*E_0\{(1-\rho^*)g(X) \}+n\Delta^*E_0[\rho^* \omega(X)g(X)].
\end{eqnarray*}
Recalling that $h(x) = 1 +\rho^*\{\omega(x) -1\} = (1 -\rho^*) + \rho^*\omega(x)$, we have 
\begin{equation*}
   E\left\{\sum_{ij} g(X_{ij}) I(X_{ij}>0) \right\} =  n\Delta^* E_0\{h(X)g(X) \}.
\end{equation*}
This completes the proof. 
$\hfill \square$

\bigskip

\cite{yuan2020zeros} define a general parameter vector $\bgamma$ of length $p$:
\begin{equation}
 \label{App_zero_gamma}
    \bgamma= \int_0^\infty \bu(x;\bnu,\btheta) dG_0(x),
\end{equation}
where $\bu(x;\bnu,\btheta)=\left(u_1(x;\bnu,\btheta),\ldots,u_p(x;\bnu,\btheta)\right)^\top$ is a given $p\times 1$ dimensional function. 
The MELE of $\bgamma$ is given by 
\begin{equation}
 \label{App_hatgamma}
    \hat{\bgamma} = \sum_{i=0}^1\sum_{j=1}^{n_{i1}} \hat p_{ij}\bu(X_{ij};\hat\bnu,\htheta).
\end{equation}

The following lemmas provide the approximation of the MELE $\hat\bgamma$ and its asymptotic property.
These lemmas help to develop the asymptotic property of the MELEs of the Gini indices. 


\begin{lemma}
 \label{App_lem3}
Assume that Conditions C1--C3 are satisfied and the true value $\nu_i^* \in (0,1)$ for $i = 0,1$.  
Let $\bgamma^*$ be the true value of $\bgamma$ and $\etab^*=(\bnu^{*\top},\rho^*,\btheta^{*\top})^\top$. Then 
\begin{eqnarray}
 \label{App_lem3_approx}
   \hat{\bgamma} &=& \frac{1}{n\Delta^*}\sum_{ij}  \frac{\bu(x;\bnu^*,\btheta^*)}{h(X_{ij})}I(X_{ij}>0) + \C(\hetab - \etab^*) + o_p(n^{-1/2}),
\end{eqnarray}
where $\C =(\C_{\bsnu}, C_{\rho}, \C_{\bstheta})$ with 
\begin{eqnarray*}
\C_{\bsnu}
&=& E_0\left\{ \frac{\partial \bu(X;\bnu^*,\btheta^*)}{\partial\bnu} \right\} + \left(\frac{w_0\bgamma^*}{\Delta^*}, \frac{w_1\bgamma^*}{\Delta^*}\right),\\
\C_\rho
&=& -E_0\left\{ \frac{\bu(X;\bnu^*,\btheta^*) \{\omega(X)-1\}}{h(X)}\right\},\\
\C_{\bstheta}
&=&
E_0\left\{\frac{\partial \bu(X;\bnu^*,\btheta^*)}{\partial\btheta}\right\}- E_0\left\{
h_1(X)\bu(X;\bnu^*,\btheta^*)\Q(X)^\top\right\}.
\end{eqnarray*}
\end{lemma}

\begin{lemma}
 \label{App_lemma_gamma}
Under the conditions of Lemma  \ref{App_lem3}, as $n\to\infty$, 
$
\sqrt{n}(\hat\bgamma- \bgamma^*) \to N({\bf 0}, \bGamma)
$
in distribution, where
\begin{eqnarray*}
\bGamma&=& \frac{1}{\Delta^*}E_0\left\{\frac{\bu(X;\bnu^*,\btheta^*)\bu(X;\bnu^*,\btheta^*)^\top}{h(X)}\right\}
- \frac{\bgamma^*\bgamma^{*\top}}{\Delta^*}\\
&&
+ \mathcal{M}_1  \A_{\bsnu}^{-1} \mathcal{M}_1 ^\top
- \frac{\mathcal{M}_2  \mathcal{M}_2^\top}{\Delta^*\rho^*(1-\rho^*)}
+ \mathcal{M}_3  \A_{\bstheta}^{-1} \mathcal{M}_3 ^\top,
\end{eqnarray*}
with
$ \A_{\bsnu} =
\diag\left\{ \frac{w_0}{\nu_0^*(1-\nu_0^*)}, \frac{w_1}{\nu_1^*(1-\nu_1^*)} \right\},~
\A_{\bstheta} =
\Delta^*(1-\rho^*)E_0\left\{h_1(X) \Q(X)\Q(X)^\top \right\}$, and
\begin{eqnarray*}
\mathcal{M}_1&=& E_0\left\{ \frac{\partial \bu(X;\bnu^*,\btheta^*)}{\partial\bnu} \right\}, \\
\mathcal{M}_2&=& E_0\left[\left\{\partial\bu(X;\bnu^*,\btheta^*)/\partial\btheta\right\}\e\right]- \rho^*\bgamma^*,\\
\mathcal{M}_3&=&   E_0\left\{\partial\bu(X;\bnu^*,\btheta^*)/\partial\btheta -
h_1(X)\bu(X;\bnu^*,\btheta^*)\Q(X)^\top\right\}.
\end{eqnarray*}
\end{lemma}

\proof
The proofs of Lemmas  \ref{App_lem3}  and  \ref{App_lemma_gamma} can be found in the supplementary material of  \cite{yuan2020zeros}. 
$\hfill \square$

\section{Proof of Theorem 1}

\subsection{Approximations of $\hat\psi_0$ and $\hat\psi_1$ \label{App_gini.supp.sec2.1}}

To develop the  asymptotic properties of $(\hat\G_0,\hat\G_1)$, 
we  first find the linear approximations of $\hat \psi_0$ and $\hat \psi_1$.  We start with $\hat \psi_0$.

Recall that 
\begin{eqnarray*}
\hat{p}_{ij} &=& \frac{1}{n_{01} +n_{11}}\left\{1 + \hat{\rho}[\exp\{ \hat\btheta^{\top} \Q(X_{ij})\}-1] \right\}^{-1}\\
&=&\frac{1}{nw_0(1-\hat\nu_0)+nw_1(1-\hat\nu_1)}\left\{1 + \hat{\rho}[\exp\{ \hat\btheta^{\top} \Q(X_{ij})\}-1] \right\}^{-1}.
\end{eqnarray*}
The MELE $\hat\psi_0$ is then given by
\begin{eqnarray*}
    \hat\psi_0 
    &=& \{nw_0(1-\hat\nu_0)+nw_1(1-\hat\nu_1)\}^{-2} \times\\
    && \sum_{ij}\sum_{ls}\frac{2I(X_{ls} \leq X_{ij})X_{ij}\cdot I(X_{ij}>0)I(X_{ls}>0)}{\left\{1 + \hat{\rho}[\exp\{ \hat\btheta^{\top} \Q(X_{ij})\}-1] \right\}\left\{1 + \hat{\rho}[\exp\{ \hat\btheta^{\top} \Q(X_{ls})\}-1] \right\}}.
\end{eqnarray*}
Note that $\hat\psi_0$ is a function of $\etab$, and hence we define 
\begin{eqnarray*}
    \psi_0(\etab)
    &=& \{nw_0(1-\nu_0)+nw_1(1- \nu_1)\}^{-2} \times\\
    && \sum_{ij}\sum_{ls}\frac{2I(X_{ls} \leq X_{ij})X_{ij}\cdot I(X_{ij}>0)I(X_{ls}>0)}{\left\{1 +  {\rho}[\exp\{  \btheta^{\top} \Q(X_{ij})\}-1] \right\}\left\{1 +  {\rho}[\exp\{  \btheta^{\top} \Q(X_{ls})\}-1] \right\}}.
\end{eqnarray*}
We then have $\hat\psi_0 = \psi_0(\hetab)$. With the definition of $\Delta^*$ and $h(x)$, we have  
\begin{eqnarray*}
\psi_0(\etab^*)= \frac{1}{(n\Delta^*)^{2}}\sum_{ij}\sum_{ls}\frac{2I(X_{ls} \leq X_{ij})X_{ij}}{h(X_{ij})h(X_{ls})}I(X_{ij}>0)I(X_{ls}>0)
\end{eqnarray*}
and $E_0\{\psi_0(\etab^*)\} = \psi_0$.

By Theorem 1 of \cite{yuan2020zeros}, we have $\hetab=\etab^*+O_p(n^{-1/2})$.
Applying the first-order Taylor expansion gives
\begin{equation}
 \label{App_psi0_taylor}
    \psi_0(\hetab)  = \psi_0(\etab^*) + \lb \frac{\partial \psi_0(\etab^*)}{\partial \etab} \rb^\top (\hetab - \etab^*) + o_p(n^{-1/2}).
\end{equation}
Define 
\begin{eqnarray*}
U(a,b) &=& \frac{I(b \leq a)a}{h(a)h(b)}I(a>0)I(b>0),~~~V_{nil} = \frac{1}{n_i}\frac{1}{n_l}\sum_{j=1}^{n_i}\sum_{s =1}^{n_l} U(X_{ij},X_{ls}),\\
V_{ni} &=& \frac{1}{n_i^2} \sum_{j=1}^{n_i}\sum_{s =1}^{n_i}2U(X_{ij},X_{is})= \frac{1}{n_i^2} \sum_{j=1}^{n_i}\sum_{s =1}^{n_i}\lb U(X_{ij},X_{is}) + U(X_{is},X_{ij})\rb,
\end{eqnarray*}
for $i,l\in \{0,1\}$ and $i \neq l$.
We then rewrite $\psi_0(\etab^*)$ as 
\begin{eqnarray*}
\psi_0(\etab^*) &=& \frac{1}{(\Delta^*)^{2}}\sum_{i=0}^1\sum_{l=0}^1 w_iw_l \frac{1}{n_i}\frac{1}{n_l}\sum_{j=1}^{n_i}\sum_{s =1}^{n_l}2U(X_{ij},X_{ls})\\
&=& \frac{1}{(\Delta^*)^{2}} \lb \sum_{i=0}^1 w_i^2 V_{ni} + \sum_{i=0}^1\sum_{l\neq i}w_iw_l 2V_{nil}\rb.
\end{eqnarray*}

Note that $V_{ni}$ is a von Mises statistic \citep{mises1947asymptotic}. 
We denote the associated U-statistic by 
$$U_{ni} = \binom{n_i}{2}^{-1}\sum_{1\leq j<s \leq n_i}\lb U(X_{ij},X_{is}) + U(X_{is},X_{ij})\rb.$$
According to \cite{Serfling1980}, the projection of $U_{ni}$ is defined as 
\begin{equation*}
    \hat U_{ni} = E\{U_{i}(X_{i1})\} + \frac{2}{n_i}\sum_{j=1}^{n_i}\Lm U_{i}(X_{ij})- E\{U_{i}(X_{i1})\} \Rm,
\end{equation*}
where $U_{i}(a) = E\{U(a,X_{i1})+U(X_{i1},a)\}$. 
It follows from \citet[p.~190 ~\& ~p.~206]{Serfling1980} that under Condition C4, 
$$
\sqrt{n_i}(\hat U_{ni} - U_{ni}) = o_p(1)~~
\mbox{and}
\sqrt{n_i}(V_{ni} - U_{ni}) = o_p(1).
$$
This leads to 
\begin{equation*}
    V_{ni} = E\{U_{i}(X_{i1})\} + \frac{2}{n_i}\sum_{j=1}^{n_i}\Lm U_{i}(X_{ij})- E\{U_{i}(X_{i1})\} \Rm + o_p(n^{-1/2}). 
\end{equation*}

When $l \neq i$, $V_{nil}$ is a two-sample U-statistic.  Define $U_{il} = E\{U(X_{i1},X_{l1})\}$, $U_{il10}(a) = E\{U(a,X_{l1})\} - U_{il}$, and $U_{il01}(a) = E\{U(X_{i1},a)\} - U_{il}$. From Theorem 12.6 in \cite{van2000asymptotic}, we have 
\begin{equation*}
    V_{nil} = U_{il} + \frac{1}{n_i}\sum_{j=1}^{n_i}U_{il10}(X_{ij}) + \frac{1}{n_l}\sum_{s=1}^{n_l} U_{il01}(X_{ls}) + o_p(n^{-1/2}).
\end{equation*}
Since 
\begin{eqnarray*}
    E\{U_{i}(X_{i1})\} = 2E\{U(X_{i1},X_{i1})\} = 2U_{ii}\\
    U_{i}(a) - E\{U_{i}(X_{i1})\} = U_{ii10}(a) + U_{ii01}(a),
\end{eqnarray*}
we have 
\begin{equation*}
    V_{ni} = 2\lb U_{ii} + \frac{1}{n_i}\sum_{j=1}^{n_i}U_{ii10}(X_{ij}) + \frac{1}{n_i}\sum_{j=1}^{n_i} U_{ii01}(X_{ij}) \rb + o_p(n^{-1/2}).
\end{equation*}
Hence, 
\begin{eqnarray}
\nonumber
  \psi_0(\etab^*) &=&
  \frac{2}{(\Delta^*)^2}\sum_{i=0}^1\sum_{l=0}^1w_iw_lU_{il}  + \frac{2}{(\Delta^*)^2}\sum_{i=0}^1\sum_{l=0}^1w_iw_l\frac{1}{n_i}\sum_{j=1}^{n_i}U_{il10}(X_{ij}) \\
  \nonumber
  && + \frac{2}{(\Delta^*)^2}\sum_{i=0}^1\sum_{l=0}^1w_iw_l\frac{1}{n_l}\sum_{s=1}^{n_l}U_{il01}(X_{ls})+o_p(n^{-1/2})\\
  \nonumber
  &=& \frac{2}{(n\Delta^*)^2}\sum_{i=0}^1\sum_{l=0}^1n_in_lU_{il}  + \frac{2}{(n\Delta^*)^2}\sum_{i=0}^1\sum_{l=0}^1 n_l\sum_{j=1}^{n_i}U_{il10}(X_{ij}) \\
   \label{App_psi0(etab*)_approx}
  && + \frac{2}{(n\Delta^*)^2}\sum_{i=0}^1\sum_{l=0}^1n_i\sum_{s=1}^{n_l}U_{il01}(X_{ls})+o_p(n^{-1/2}).
\end{eqnarray}

We now simplify each term in  \eqref{App_psi0(etab*)_approx}. 
With Lemma  \ref{App_lem1} and  the definition of $U_{il}$, we have 
\begin{eqnarray*}
  \sum_{l=0}^1n_lE\{U(X_{i1},X_{l1})|X_{i1}\} &=&  \sum_{l=0}^1n_lE\left\{\frac{I(X_{l1}\leq X_{i1})X_{i1}}{h(X_{i1})h(X_{l1})}I(X_{i1}>0)I(X_{l1}>0)|X_{i1}\right\}\\
  &=& n\Delta^*\frac{X_{i1}}{h(X_{i1})}I(X_{i1}>0) E_0\left\{I(X\leq X_{i1})\right\}\\
  &=& n\Delta^*\frac{X_{i1}G_0(X_{i1})}{h(X_{i1})}I(X_{i1}>0) .
\end{eqnarray*}
Hence, 
\begin{eqnarray*}
   \frac{2}{(n\Delta^*)^2}\sum_{i=0}^1\sum_{l=0}^1n_in_lU_{il}
   &=& \frac{2}{(n\Delta^*)^2}\sum_{i=0}^1n_iE\Lm \sum_{l=0}^1n_lE\{U(X_{i1},X_{l1})|X_i\}\Rm\\
   &=& \frac{2}{n\Delta^*}\sum_{i=0}^1n_i E\lb \frac{X_{i1} G_0(X_{i1})}{h(X_{i1})}I(X_{i1}>0)\rb.
\end{eqnarray*}
Using the result in Lemma  \ref{App_lem1}, we have 
\begin{equation*}
     \frac{2}{(n\Delta^*)^2}\sum_{i=0}^1\sum_{l=0}^1n_in_lU_{il} = 2E_0\{XG_0(X)\} =\psi_0.
\end{equation*}

We move to the second term of $\psi_0(\etab^*)$ in  \eqref{App_psi0(etab*)_approx}. 
Recall that $$U_{il10}(a) = E\{U(a,X_{l1})\} - E\{U(X_{i1},X_{l1})\}.$$
We then have 
\begin{eqnarray*}
    \sum_{l=0}^1 n_l U_{il10}(X_{ij})
    &=& \sum_{l=0}^1 n_l \left\{E\{U(X_{ij},X_{l1})|X_{ij}\} - E[E\{U(X_{ij},X_{l1})|X_{ij}\}]\right\}\\
    &=& n\Delta^*\left[\frac{X_{ij}G_0(X_{ij})}{h(X_{ij})}I(X_{ij}>0) - E\left\{ \frac{X_{ij}G_0(X_{ij})}{h(X_{ij})}I(X_{ij}>0)\right\}\right].
\end{eqnarray*}
This leads to 
\begin{eqnarray*}
    &&\frac{2}{(n\Delta^*)^2}\sum_{i=0}^1\sum_{l=0}^1 n_l\sum_{j=1}^{n_i}U_{il10}(X_{ij})\\
    &=& 
    \frac{2}{n\Delta^*}\sum_{ij} \left[\frac{X_{ij}G_0(X_{ij})}{h(X_{ij})}I(X_{ij}>0)  - E\left\{ \frac{X_{ij}G_0(X_{ij})}{h(X_{ij})}I(X_{ij}>0)\right\}\right]\\
    &=&\frac{2}{n\Delta^*}\sum_{ij} \frac{I(X_{ij}>0)}{h(X_{ij})} X_{ij}G_0(X_{ij}) -  \psi_0.
\end{eqnarray*}
Similarly, with the definition of $U_{il01}(a)$, we have
\begin{eqnarray*}
    \sum_{i=0}^1 n_i U_{il01}(X_{ls})&=& \sum_{i=0}^1 n_i \left\{E\{U(X_{i1},X_{ls})|X_{ls}\} - E[E\{U(X_{i1},X_{ls})|X_{ls}\}]\right\}.
\end{eqnarray*}
Note that 
\begin{eqnarray*}
  \sum_{i=0}^1 n_i E\{U(X_{i1},X_{ls})|X_{ls}\} &=&  \sum_{i=0}^1 n_iE\left\{\frac{I(X_{ls}\leq X_{i1})X_{i1}}{h(X_{i1})h(X_{ls})}I(X_{i1}>0)I(X_{ls}>0)|X_{ls}\right\}\\
  &=& n\Delta^*\frac{I(X_{ls}>0)}{h(X_{ls})} E\left\{XI(X_{ls} \leq X)\right\}\\
  &=& n\Delta^*\frac{I(X_{ls}>0)}{h(X_{ls})}\int_{X_{ls}}^{\infty} x dG_0(x).
\end{eqnarray*}
Together with the result of Lemma  \ref{App_lem1}, we have 
\begin{eqnarray*}
    &&\frac{2}{(n\Delta^*)^2}\sum_{i=0}^1 \sum_{l =0}^1 n_i \sum_{s = 1}^{n_l} U_{il01}(X_{ls}) \\
    &=& 
    \frac{2}{n\Delta^*}\sum_{ls}\left[ \frac{I(X_{ls}>0)}{h(X_{ls})}\int_{X_{ls}}^{\infty} xdG_0(x) - E\left\{ \frac{I(X_{ls}>0)}{h(X_{ls})}\int_{X_{ls}}^{\infty} x dG_0(x)\right\}\right]\\
    &=& \frac{2}{n\Delta^*}\sum_{ls} \frac{I(X_{ls}>0)}{h(X_{ls})}\int_{X_{ls}}^{\infty} xdG_0(x) - \psi_0.
\end{eqnarray*}
For $a > 0$, we define  the function
\begin{equation*}
    H_0(a) = \left\{aG_0(a)+ \int_{a}^{\infty} xdG_0(x)\right\}.
\end{equation*}
The approximation of $\psi_0(\etab^*)$ is then given by 
\begin{eqnarray*}
\psi_0(\etab^*) = \frac{1}{n\Delta^*}\sum_{ij} \frac{I(X_{ij}>0)}{h(X_{ij})}\cdot \{2H_0(X_{ij})\} -  \psi_0 + o_p(n^{-1/2}).
\end{eqnarray*}

We also need the first derivative of $\psi_0(\etab)$ when finding the approximation of $\psi_0(\hetab)$.
We take the first derivative of $\psi_0(\etab)$ with respect to $\etab$ and evaluate the derivative at the true value $\etab^*$. This leads to 
\begin{eqnarray*}
\frac{\partial \psi_0(\etab^*)}{\partial \bnu}
&=&
\frac{2}{\Delta^*} \psi_0(\etab^*)\ba{c}w_0 \\
  w_1 \ea,\\
\frac{\partial \psi_0(\etab^*)}{\partial \rho}
&=&
-\frac{2}{(n\Delta^*)^2}\sum_{ij}\sum_{ls}\lb \frac{\omega(X_{ij})-1}{h(X_{ij})^2h(X_{ls})}+\frac{\omega(X_{ls})-1}{h(X_{ij})h(X_{ls})^2}\rb \\
&&\times I(X_{ls} \leq X_{ij})X_{ij}I(X_{ij}>0)I(X_{ls} >0),\\
\frac{\partial \psi_0(\etab^*)}{\partial \btheta}
&=&
-\frac{2}{(n\Delta^*)^2}\sum_{ij}\sum_{ls}\lb \frac{\rho^*\omega(X_{ij})\Q(X_{ij})}{h(X_{ij})^2h(X_{ls})}+\frac{\rho^*\omega(X_{ls})\Q(X_{ls})}{h(X_{ij})h(X_{ls})^2}\rb\\
&& \times I(X_{ls} \leq X_{ij})X_{ij}I(X_{ij}>0)I(X_{ls} >0),\\
&=& -\frac{2}{(n\Delta^*)^2}\sum_{ij}\sum_{ls}\lb h_1(X_{ij})\Q(X_{ij})+h_1(X_{ls})\Q(X_{ls})\rb\\
&& \times \frac{I(X_{ls} \leq X_{ij})X_{ij}I(X_{ij}>0)I(X_{ls} >0)}{h(X_{ij})h(X_{ls})}.
\end{eqnarray*}
By the law of large numbers, we have 
\begin{equation*}
     \frac{\partial \psi_0(\etab^*)}{\partial \etab} = E\lb \frac{\partial \psi_0(\etab^*)}{\partial \etab} \rb + o_p(1) = \C_0 + o_p(1),
\end{equation*}
with $\C_0 = (\C_{0\bsnu}^\top,\C_{0\rho},\C_{0\bstheta}^\top)^\top$.

Since  $E\{\psi_0(\etab^*)\} =  \psi_0$, we have 
\begin{equation*}
    \C_{0\bsnu} = \frac{2 \psi_0}{\Delta^*}(w_0,w_1)^\top.
\end{equation*}
For $\C_{0\rho}$, 
\begin{eqnarray*}
\C_{0\rho}&=& E\lb \frac{\partial \psi_0(\etab^*)}{\partial \rho} \rb \\
&=&-\frac{2}{(n\Delta^*)^2} E\left[\sum_{ij} E\lb \sum_{ls}\frac{\omega(X_{ij})-1}{h(X_{ij})^2h(X_{ls})} I(X_{ls} \leq X_{ij})X_{ij}I(X_{ij}>0)I(X_{ls} >0)|X_{ij} \rb \right]\\ 
&&-\frac{2}{(n\Delta^*)^2} E\left[\sum_{ls} E\lb \sum_{ij}\frac{\omega(X_{ls})-1}{h(X_{ij})h(X_{ls})^2} I(X_{ls} \leq X_{ij})X_{ij}I(X_{ij}>0)I(X_{ls} >0)|X_{ls} \rb \right]\\
&=&-\frac{2}{n\Delta^*} E\left[\sum_{ij} \frac{\omega(X_{ij})-1}{h(X_{ij})^2} X_{ij}I(X_{ij}>0)G_0(x)\right]\\ 
&&-\frac{2}{n\Delta^*} E\left[\sum_{ls} \frac{\omega(X_{ls})-1}{h(X_{ls})^2}I(X_{ls} >0)E_0\lb  I(X_{ls} \leq X)X\rb \right]\\
&=& -\frac{2}{n\Delta^*} E\left[\sum_{ij} \frac{\omega(X_{ij})-1}{h(X_{ij})^2}I(X_{ij}>0)H_0(X_{ij})\right]\\ 
&=& -2E_0\left\{ \frac{H_0(X) \{\omega(X)-1\}}{h(X)}\right\}.
\end{eqnarray*}
The expression for $\C_{0\bstheta}$ can be found in a similar manner: 
\begin{equation*}
    \C_{0\bstheta} = -2 E_0\left\{
h_1(X)H_0(X)\Q(X)\right\}.
\end{equation*}
The details  are omitted here. 

It can be verified that the matrix $\C_0^\top$ is the same as the matrix $\C$ in  \eqref{App_lem3_approx} when we set $\bu(x;\bnu,\btheta) = 2H_0(x)$ in  \eqref{App_zero_gamma}. 
According to Lemma  \ref{App_lem3}, the expression in  \eqref{App_psi0_taylor} can be further written as 
\begin{eqnarray}
\nonumber
   \hat\psi_0 &=& \frac{1}{n\Delta^*}\sum_{ij} \frac{I(X_{ij}>0)}{h(X_{ij})}\cdot \{2H_0(X_{ij})\} -  \psi_0
   + \C_0^\top(\hetab - \etab^*)+o_p(n^{-1/2})\\
    \label{App_hatpsi0_approx}
   &=&\sum_{i=0}^1\sum_{j = 1}^{n_{i1}}\hat{p}_{ij}\{2H_0(X_{ij})\}- \psi_0 + o_p(n^{-1/2}).
\end{eqnarray}
The remaining term $o_p(n^{-1/2})$ is introduced by the projection of the von Mises statistic and the U-statistic when approximating $\psi_0(\etab^*)$. 

Define 
\begin{equation*}
    \mathcal{H}_0(a) = 2H_0(a) -  \psi_0.
\end{equation*}
With the natural constraint $ \sum_{i=0}^1\sum_{j = 1}^{n_{i1}} p_{ij} =1$,
Equation  \eqref{App_hatpsi0_approx} implies 
\begin{equation}
 \label{App_approx.psi0}
    \hat\psi_0 = \sum_{i=0}^1\sum_{j = 1}^{n_{i1}}\hat{p}_{ij}\mathcal{H}_0(X_{ij}) + o_p(n^{-1/2}).
\end{equation}


Next, we consider the approximation of MELE $\hat\psi_1$. Recall that 
\begin{equation*}
    \hat\psi_1 = \sum_{i =0}^1\sum_{j =1}^{n_{i1}} \hat{p}_{ij}\omega(X_{ij};\htheta) X_{ij}\lb 2 \sum_{l =0}^1 \sum_{s = 1}^{n_{l1}} \hat{p}_{ls}\omega(X_{ls};\htheta)I(X_{ls} \leq X_{ij}) \rb.
\end{equation*} 
With the definition of $\hat p_{ij}$, the MELE $\hat\psi_i$ can be written as 
\begin{eqnarray}
\nonumber
    \hat\psi_1
    &=& \{nw_0(1-\hat\nu_0)+nw_1(1-\hat\nu_1)\}^{-2}\\
     \label{App_hatpsi12}
    && \times \sum_{ij}\sum_{ls}\frac{2I(X_{ls} \leq X_{ij})X_{ij}\omega(X_{ij};\htheta)\omega(X_{ls};\htheta)I(X_{ij}>0)I(X_{ls}>0)}{\left\{1 + \hat{\rho}[\exp\{ \hat\btheta^{\top} \Q(X_{ij})\}-1] \right\}\left\{1 + \hat{\rho}[\exp\{ \hat\btheta^{\top} \Q(X_{ls})\}-1] \right\}}.
\end{eqnarray}

Define 
\begin{eqnarray*}
\tilde{U}(a,b) &=& \frac{I(b \leq a)a}{h(a)h(b)}\omega(a)\omega(b)I(a>0)I(b>0),~~~\tilde{V}_{nil} = \frac{1}{n_i}\frac{1}{n_l}\sum_{j=1}^{n_i}\sum_{s =1}^{n_l} \tilde{U}(X_{ij},X_{ls}),\\
\tilde{V}_{ni} &=& \frac{1}{n_i^2} \sum_{j=1}^{n_i}\sum_{s =1}^{n_i}2\tilde{U}(X_{ij},X_{is})= \frac{1}{n_i^2} \sum_{j=1}^{n_i}\sum_{s =1}^{n_i}\lb \tilde{U}(X_{ij},X_{is}) + \tilde{U}(X_{is},X_{ij})\rb,
\end{eqnarray*}
for $i,l \in\{0,1\}$.
We use $\psi_1(\hat\etab)$ to denote $\hat\psi_1$ and have 
\begin{eqnarray*}
\psi_1(\etab^*)&=& \frac{1}{(n\Delta^*)^2}\sum_{ij}\sum_{ls}2\tilde{U}(X_{ij},X_{ls})\\
&=& \frac{1}{(\Delta^*)^{2}}\sum_{i=0}^1\sum_{l=0}^1 w_iw_l \frac{1}{n_i}\frac{1}{n_l}\sum_{j=1}^{n_i}\sum_{s =1}^{n_l}2\tilde{U}(X_{ij},X_{ls})\\
&=& \frac{1}{(\Delta^*)^{2}} \lb \sum_{i=0}^1 w_i^2 \tilde{V}_{ni} + \sum_{i=0}^1\sum_{l\neq i}w_iw_l 2\tilde{V}_{nil}\rb.
\end{eqnarray*}

Note that  $\tilde{V}_{ni}$ is a von Mises statistic and $\tilde{V}_{nil}$ is a two-sample U-statistic. 
Using the technique used to obtain the approximation of $\psi_0(\etab^*)$ in  \eqref{App_psi0(etab*)_approx}, we have  
\begin{eqnarray*}
  \psi_1(\etab^*)
  &=& \frac{2}{(n\Delta^*)^2}\sum_{i=0}^1\sum_{l=0}^1n_in_l\tilde{U}_{il}  + \frac{2}{(n\Delta^*)^2}\sum_{i=0}^1\sum_{l=0}^1 n_l\sum_{j=1}^{n_i}\tilde{U}_{il10}(X_{ij}) \\
  && + \frac{2}{(n\Delta^*)^2}\sum_{i=0}^1\sum_{l=0}^1n_i\sum_{s=1}^{n_l}\tilde{U}_{il01}(X_{ls})+o_p(n^{-1/2}),
\end{eqnarray*}
where $\tilde{U}_{il} = E\{\tilde{U}(X_{i1},X_{l1})\}$, $\tilde{U}_{il10}(a) = E\{\tilde{U}(a,X_{l1})\} - \tilde{U}_{il}$, and $\tilde{U}_{il01}(a) = E\{\tilde{U}(X_{i1},a)\} - \tilde{U}_{il}$.

With Lemma  \ref{App_lem1} and the definition of $\tilde{U}(a,b)$, we have 
\begin{eqnarray*}
 &&\sum_{l=0}^1n_lE\{\tilde{U}(X_{ij},X_{l1})|X_{ij}\}\\
 &=&  \sum_{l=0}^1n_lE\left\{\frac{I(X_{l1}\leq X_{ij})X_{ij}}{h(X_{ij})h(X_{l1})}\omega(X_{ij})\omega(X_{l1})I(X_{ij}>0)I(X_{l1}>0)|X_{ij}\right\}\\
  &=& n\Delta^*\frac{X_{ij}\omega(X_{ij})}{h(X_{ij})}I(X_{ij}>0) E_0\left\{\omega(X)I(X\leq X_{ij})\right\}\\
  &=& n\Delta^*\frac{X_{ij}\omega(X_{ij})G_1(X_{ij})}{h(X_{ij})}I(X_{ij}>0)
\end{eqnarray*}
and
\begin{eqnarray*}
  &&\sum_{i=0}^1 n_i E\{\tilde{U}(X_{i1},X_{ls})|X_{ls}\}\\
  &=&  \sum_{i=0}^1 n_iE\left\{\frac{I(X_{ls}\leq X_{i1})X_{i1}}{h(X_{i1})h(X_{ls})}\omega(X_{i1})\omega(X_{ls})I(X_{i1}>0)I(X_{ls}>0)|X_{ls}\right\}\\
  &=& n\Delta^*\frac{\omega(X_{ls})I(X_{ls}>0)}{h(X_{ls})} E_0\left\{X\omega(X)I(X_{ls} \leq X)\right\}\\
  &=& n\Delta^*\frac{\omega(X_{ls})I(X_{ls}>0)}{h(X_{ls})}\int_{X_{ls}}^{\infty} x dG_1(x).
\end{eqnarray*}
It follows that
\begin{eqnarray*}
\frac{2}{(n\Delta^*)^2}\sum_{i=0}^1\sum_{l=0}^1n_in_l\tilde{U}_{il}
&=&\frac{2}{n\Delta^*}\sum_{i=0}^1 n_iE\left\{\frac{X_{i1}\omega(X_{i1})G_1(X_{i1})}{h(X_{i1})}I(X_{i1}>0) \right\}\\ 
&=& 2E_0\{ X\omega(X)G_1(X)\}\\
&=&  \psi_1,\\
\frac{2}{(n\Delta^*)^2}\sum_{i=0}^1\sum_{l=0}^1 n_l\sum_{j=1}^{n_i}\tilde{U}_{il10}(X_{ij}) 
&=& \frac{2}{n\Delta^*}\sum_{ij}\frac{\omega(X_{ij})I(X_{ij}>0)}{h(X_{ij})}X_{ij}G_1(X_{ij})  -  \psi_1,\\
\frac{2}{(n\Delta^*)^2}\sum_{i=0}^1\sum_{l=0}^1n_i\sum_{s=1}^{n_l}\tilde{U}_{il01}(X_{ls})
&=&\frac{2}{n\Delta^*}\sum_{ls}\frac{\omega(X_{ls})I(X_{ls}>0)}{h(X_{ls})}\int_{X_{ls}}^{\infty} x dG_1(x) - \psi_1.
\end{eqnarray*}

Hence, $\psi_1(\etab^*)$ is given by 
\begin{equation*}
   \psi_1(\etab^*) = \frac{1}{n\Delta^*}\sum_{ij} \frac{\omega(X_{ij}) I(X_{ij}>0)}{h(X_{ij})}\cdot \{2H_1(X_{ij})\} -  \psi_1 + o_p(n^{-1/2}),
\end{equation*}
where $H_1(a) = aG_1(a)+ \int_{a}^{\infty} xdG_1(x)$ for $a>0$.

Applying the first-order Taylor expansion to $\hat\psi_1$ in  \eqref{App_hatpsi12} yields
\begin{eqnarray*}
\psi_1(\hetab) = \psi_1(\etab^*) + \lb \frac{\partial \psi_0(\etab^*)}{\partial \etab} \rb^\top (\hetab - \etab^*) + o_p(n^{-1/2}).
\end{eqnarray*}
With the law of large numbers, we have 
\begin{equation*}
     \frac{\partial \psi_1(\etab^*)}{\partial \etab} = E\lb \frac{\partial \psi_1(\etab^*)}{\partial \etab} \rb + o_p(1) = \C_1 + o_p(1)
\end{equation*}
with $\C_1 = (\C_{1\bsnu}^\top,\C_{1\rho},\C_{1\bstheta}^\top)^\top$ and
\begin{equation*}
 \C_{1\bsnu} = \frac{2 \psi_1}{\Delta^*}
  \ba{c}w_0 \\
  w_1 \ea,  \C_{1\rho} =
  -2E_0\left\{ \frac{H_1(X) \{\omega(X)-1\}}{h(X)}\right\},\C_{1\bstheta} = 2 E_0\left\{
h_0(X)H_1(X)\Q(X)\right\}.
\end{equation*}
The expression of each element in $\C_1$ can be found similarly to the derivation of $\C_0$; we omit the details.
By setting $\bu(x;\bnu,\btheta) = 2\omega(x;\btheta)H_1(x)$ in  \eqref{App_zero_gamma}, we can verify that the matrix $\C_1^\top$ is the same as the matrix $\C$ in  \eqref{App_lem3_approx}.
Hence, the approximation is given by 
\begin{equation*}
    \psi_1(\hetab) = \frac{1}{n\Delta^*}\sum_{ij} \frac{\omega(X_{ij}) I(X_{ij}>0)}{h(X_{ij})}\cdot\{2H_1(X_{ij})\} -  \psi_1 + \C_1^\top (\hetab - \etab^*) + o_p(n^{-1/2}).
\end{equation*}

With Lemma  \ref{App_lem3} and the natural constraint $\sum_{i=0}^1\sum_{j=1}^{n_i}\hat p_{ij}\omega(X_{ij};\htheta) = 1$, the above approximation equation implies 
\begin{equation}
 \label{App_approx.psi1}
    \hat\psi_1 = \sum_{i=0}^1\sum_{j = 1}^{n_{i1}}\hat{p}_{ij}\omega(X_{ij};\htheta)\mathcal{H}_1(X_{ij}) + o_p(n^{-1/2}),
\end{equation}
where we define $\mathcal{H}_1(a) = 2H_1(a) - \psi_1$ for $a>0$.

\subsection{Asymptotic properties  of  $\hat \G_0$ and $\hat \G_1$}
In this section, we use the approximations of $\hat\psi_0$ and $\hat\psi_1$ developed in Section  \ref{App_gini.supp.sec2.1} 
and the results in Lemmas  \ref{App_lem3} and  \ref{App_lemma_gamma}
to derive the asymptotic properties of $\hat \G_0$ and $\hat \G_1$. 

Recall that 
$$
\hat\G_0 =(2\hat\nu_0 - 1)+(1-\hat{\nu}_0) \frac{\hat\psi_0}{\hat m_0}, 
$$
where $\hat m_0 = \sum_{i=0}^1\sum_{j=1}^{n_{i1}} \hat{p}_{ij} X_{ij}$ and the approximation of $\hat\psi_0$ is in  \eqref{App_approx.psi0}. 
Then 
\begin{eqnarray}
\nonumber
\hat\G_0 &=& \frac{(2\hat\nu_0 - 1)\hat m_0 + (1-\hat\nu_0)\hat\psi_0}{\hat m_0}\\
&=& \frac{  \sum_{i=0}^1\sum_{j = 1}^{n_{i1}}\hat{p}_{ij}\left\{ (2\hat\nu_0 - 1)X_{ij} + (1-\hat\nu_0)\mathcal{H}_0(X_{ij}) \right\} }{  \sum_{i=0}^1\sum_{j = 1}^{n_{i1}}\hat{p}_{ij}X_{ij}}+ o_p(n^{-1/2}). 
 \label{App_hatG0.approx1}
\end{eqnarray}
Similarly, 
\begin{eqnarray}
\hat\G_1 &=&\frac{ \sum_{i=0}^1\sum_{j = 1}^{n_{i1}}\hat{p}_{ij}\left\{ (2\hat\nu_1 - 1)\omega(X_{ij};\htheta)X_{ij} + (1-\hat\nu_1)\mathcal{H}_1(X_{ij};\htheta) \right\}  }{  \sum_{i=0}^1\sum_{j = 1}^{n_{i1}}\hat{p}_{ij}\omega(X_{ij};\htheta)X_{ij} }+ o_p(n^{-1/2}).
 \label{App_hatG1.approx1}
\end{eqnarray}

Note that the numerators and denominators of the leading terms in  \eqref{App_hatG0.approx1} and  \eqref{App_hatG1.approx1}
all have the forms in  \eqref{App_hatgamma} with $\bu(\cdot;\cdot)$ taking some specific forms. 
We define these specific $\bu(x;\bnu,\btheta)$ as 
\begin{eqnarray}
 \label{App_Gini_u}
\bu(x;\bnu,\btheta) &=& \left(x,u_0(x;\bnu),\omega(x;\btheta)x,\omega(x;\btheta)u_1(x;\bnu)\right)^\top 
\end{eqnarray}
with 
\begin{eqnarray}
 \label{App_gini.u01}
u_0(x;\bnu)=(2\nu_0 -1)x + (1-\nu_0)\mathcal{H}_0(x)
~~
\mbox{and}
~~
u_1(x;\bnu)=(2\nu_1 -1)x + (1-\nu_1)\mathcal{H}_1(x).
\end{eqnarray}
Further, we define 
\begin{eqnarray}
 \label{App_Gini_gamma}
\bgamma =\int_0^{\infty}\bu(x;\bnu,\btheta)dG_0(x) = (\gamma_1, \gamma_2,\gamma_3,\gamma_4)^\top
\end{eqnarray}
and 
\begin{eqnarray}
 \label{App_Gini_hatgamma}
\hat\bgamma = \int_0^{\infty}\bu(x;\hat\bnu,\hat\btheta)d\hat G_0(x) = (\hat\gamma_1, \hat\gamma_2,\hat\gamma_3,\hat\gamma_4)^\top. 
\end{eqnarray}
Then we have 
\begin{eqnarray}
\hat\G_0 =\hat\gamma_1/\hat\gamma_2+ o_p(n^{-1/2})~~
 \mbox{ and } 
 ~~
 \hat\G_1 =\hat\gamma_3/\hat\gamma_4+ o_p(n^{-1/2}). 
\end{eqnarray}
Hence, the joint limiting distribution of $\sqrt{n}(\hat\G_0-\G_0^*,\hat\G_1-\G_1^*)$
is determined by that of $\sqrt{n}(\hat\bgamma-\bgamma^*)$, 
where the $\G_i^*$ are the true values of $\G_i$ for $i =0,1$, 
and  
\begin{eqnarray*}
\bgamma^* &=& \int_0^{\infty}\bu(x;\bnu^*,\btheta^*)dG_0(x) = (m_0, m_0\G_0^*, m_1, m_1\G_1^*)^\top
\end{eqnarray*}
is the true value of $\bgamma$. 

Let $$
\tilde{\bu}(x) 
= (\tilde{\bu}_0(x)^\top,\tilde{\bu}_1(x)^\top)^\top = (-\rho^*(x,u_0(x;\bnu^*)) ,(1-\rho^*)(x,u_1(x;\bnu^*)) )^\top.
$$
Applying Lemma  \ref{App_lemma_gamma}, we have,  as $n$ goes to infinity,
\begin{equation*}
    \sqrt{n}(\hat{\bgamma} - \bgamma^*) \to N({\bf 0}, \bGamma)
\end{equation*}
in distribution with
\begin{equation}
     \label{App_bGamma}
    \bGamma= \frac{1}{\Delta^*}E_0\left\{\frac{\bu(X;\bnu^*,\btheta^*)\bu(X;\bnu^*,\btheta^*)^\top}{h(X)}\right\}
- \frac{\bgamma^*\bgamma^{*\top}}{\Delta^*}
+ \mathcal{M}_1  \A_{\bsnu}^{-1} \mathcal{M}_1 ^\top
- \frac{\mathcal{M}_2  \mathcal{M}_2^\top}{\Delta^*\rho^*(1-\rho^*)}
+ \mathcal{M}_3  \A_{\bstheta}^{-1} \mathcal{M}_3 ^\top,
\end{equation}
where $ \A_{\bsnu}$ and $\A_{\bstheta}$ are provided in Lemma  \ref{App_lemma_gamma}, and 
 \begin{equation*}
\mathcal{M}_1 =
\ba{cc} 0& 0\\ 2m_0 -  \psi_0&0\\
0&0 \\ 0&2m_1 -  \psi_1\ea,~~~
\mathcal{M}_2 =\ba{c} -\rho^*m_0\\-\rho^*(m_0\G_0^*)\\ (1-\rho^*)m_1\\(1-\rho^*)(m_1\G_1^*)\ea 
= E_0\{\tilde{\bu}(X)\},
\end{equation*}
\begin{equation*}
\mathcal{M}_3= 
\ba{c}
-E_0\{h_1(X)X\Q(X)^\top\}\\
-E_0\{h_1(X)u_0(X;\bnu^*)\Q(X)^\top\}\\
E_0\{h_0(X)\omega(X)X\Q(X)^\top\}\\
E_0\{h_0(X)\omega(X)u_1(X;\bnu^*)\Q(X)^\top\} \ea\\
=\frac{1}{\rho^*}E_0\{h_1(X)\tilde{\bu}(X)\Q(X)^\top\}. 
\end{equation*}
Note that we have used the form of $\bu(x;\bnu,\btheta)$ in  \eqref{App_Gini_u} to obtain the simplified forms of $\mathcal{M}_1$, $\mathcal{M}_2$, and $\mathcal{M}_3$. 

Let $\g(\bgamma)= (\gamma_2/\gamma_1,\gamma_4/\gamma_3)^\top = (\G_0,\G_1)^\top$. 
Then $\g(\hat\bgamma) = (\hat\G_0,\hat \G_1)^\top$ and $\g(\bgamma^*) = (\G_0^*,\G_1^*)^\top$. Using the Delta method,   we have,  as $n\to\infty$, 
\begin{equation*}
    \sqrt{n}\left(\hat\G_0 - \G_0^*\right)  \to N({\bf 0},\bSigma)
\end{equation*}
in distribution, where $\bSigma= \J \bGamma \J^\top$ and 
\begin{equation}
 \label{App_gini.Jmat}
    \J = \frac{\partial \g(\bgamma^*)}{\partial \bgamma } = \ba{cc} \J_0^\top & {\bf 0}\\
    {\bf 0} & \J_1^\top\ea = 
    \ba{cccc}
    -\frac{\G_0^*}{m_0}&\frac{1}{m_0}&0&0\\
    0&0&-\frac{\G_1^*}{m_1}&\frac{1}{m_1}
    \ea. 
\end{equation}

To finish the proof of Theorem 1, we use the forms of $\bGamma$ in  \eqref{App_bGamma} and $\J$ in  \eqref{App_gini.Jmat}
to simplify $\bSigma$. 
Note that 
\begin{eqnarray*}
\J_0^\top \ba{c} m_0\\m_0\G_0^*\ea = -\frac{\G_0^*}{m_0}\cdot m_0 + \frac{1}{m_0}\cdot (m_0\G_0^*) = 0,\\
\J_1^\top \ba{c} m_1\\m_1\G_1^*\ea = -\frac{\G_1^*}{m_1}\cdot m_1 + \frac{1}{m_1}\cdot (m_1\G_1^*) = 0.
\end{eqnarray*}
This leads to 
\begin{equation}
 \label{App_result_bSigma1}
    \J \lb \frac{\bgamma^*\bgamma^{*\top}}{\Delta^*}\rb \J^\top = \J \lb- \frac{\mathcal{M}_2  \mathcal{M}_2^\top}{\Delta^*\rho^*(1-\rho^*)}\rb \J^\top = {\bf 0}.
\end{equation}

With the fact that 
\begin{eqnarray*}
\J_0^\top \ba{c} 0\\2m_0 -  \psi_0\ea = \frac{1 - \G_0^*}{1-\nu_0}~~\text{and}~~
\J_1^\top \ba{c} 0\\2m_1 -  \psi_1\ea = \frac{1 - \G_1^*}{1-\nu_1},
\end{eqnarray*}
we have 
\begin{equation*}
    \J \mathcal{M}_1 = \diag\lb \frac{1 - \G_0^*}{1-\nu_0}, \frac{1 - \G_1^*}{1-\nu_1}\rb.
\end{equation*}
Hence, 
\begin{equation}
 \label{App_result_bSigma2}
    \J (\mathcal{M}_1\A_{\bsnu}^{-1}\mathcal{M}_1) \J^\top= \diag\lb \frac{\nu_0^*(1 - \G_0^*)^2}{\Delta^*(1-\rho^*)},\frac{\nu_1^*(1 - \G_1^*)^2}{\Delta^*\rho^*}\rb.
\end{equation}

Substituting  \eqref{App_bGamma} and  \eqref{App_gini.Jmat} into $\bSigma$ and using    \eqref{App_result_bSigma1}-- \eqref{App_result_bSigma2},  $\bSigma$  
has the following simplified form: 
\begin{eqnarray}
\nonumber
\bSigma &=& \frac{1}{\Delta^*}\J \Bigg[ E_0\left\{\frac{\bu(X;\bnu^*,\btheta^*)\bu(X;\bnu^*,\btheta^*)^\top}{h(X)}\right\} + \frac{1}{(\rho^*)^2}\B \Bigg] \J^\top\\
 \label{App_final_bSigma}
&& + \diag\lb \frac{\nu_0^*(1 - \G_0^*)^2}{\Delta^*(1-\rho^*)},\frac{\nu_1^*(1 - \G_1^*)^2}{\Delta^*\rho^*}\rb,
\end{eqnarray}
where 
\begin{equation*}
    \B = E_0\{h_1(X)\tilde{\bu}(X)\Q(X)^\top\}\A_{\btheta}^{-1} E_0\{h_1(X)\Q(X)\tilde{\bu}(X)^\top\}, 
\end{equation*}
as claimed in Theorem 1. This completes the proof. 

\section{Proof of Theorem 2}

The proof of Theorem 2 is similar to that of Theorem 1. 
The  results of \cite{li2018comparison} are helpful for this proof.

\section{Proof of Theorem 3}

We start with (a). 
Recall that 
the nonparametric estimator of the Gini index for sample $i = 0,1$ is defined as 
$$
\tilde\G_i =(2\hat\nu_i-1)+(1-\hat\nu_i){\tilde\psi_i}/{\tilde m_i},
$$
where 
\begin{eqnarray*}
\tilde m_i=  {n_{i1}}^{-1}\sum_{j=1}^{n_{i1}}X_{ij} ~~
\mbox{and}
~~
\tilde \psi_i=  \int_0^{\infty} \{2x \tilde G_i(x)\}d \tilde G_i(x).
\end{eqnarray*}
After some algebra, we have 
\begin{equation}
    \tilde{\G}_i = \frac{n_i^{-1}\sum_{j = 1}^{n_i} \{2\tilde{F}_i(X_{ij}) -1\}X_{ij}}{\tilde{\mu}_i}, 
\end{equation}
where $\tilde{\mu}_i = {n_i}^{-1}\sum_{j=1}^{n_i}X_{ij}$ and $\tilde{F}_i(x) = {n_i}^{-1}\sum_{j=1}^{n_i}I(X_{ij} \leq x)$ are the sample mean and the empirical CDF based on  sample $i$.

Applying Theorem 1 of \cite{qin2010empirical}, 
we have 
\begin{equation*}
    \sqrt{n} \ba{c} \tilde{\G}_0 - \G_0^*\\ \tilde{\G}_1 - \G_1^*\ea  \to N ({\bf 0}, \bSigma_{non})
\end{equation*}
in distribution with 
\begin{equation}
 \label{App_bSigma_non}
    \bSigma_{non} =  \ba{cc} \sigma_0^2&0\\ 0& \sigma_1^2 \ea,
\end{equation} 
where 
\begin{equation}
 \label{App_gini.sigma.i1}
    \sigma_i^2 =\frac{ Var_{F_i}\{u_i(X)-\G_i^* X \} }{w_i\mu_i^{2}},
\end{equation}
where $Var_{F_i}$ means the variance is taken with respect to $F_i$ and $u_i(x)=u_i(x;\bnu^*)$ with $u_i(x;\bnu)$ defined in  \eqref{App_gini.u01}.

We now show that $   \bSigma_{non} $ has the form claimed in (a). 
Note that 
\begin{eqnarray}
 \label{App_gini.non.var1}
Var_{F_i}\{u_i(X)-\G_i^* X \}&=&E_{F_i}[\{u_i(X)-\G_i^* X \}^2]- [E_{F_i}\{u_i(X)-\G_i^* X \}]^2,
\end{eqnarray}
where $E_{F_i}$ indicates that the expectation is taken with respect to $F_i$. 
After some calculus work, we can show that 
\begin{equation}
 \label{App_gini.non.exp1}
E_{F_i}\{u_i(X)-\G_i^* X \}=\nu_i^*(1-\nu_i^*)(2m_i-\psi_i). 
\end{equation}
For 
$E_{F_i}[\{u_i(X)-\G_i^* X \}^2]$, 
we have under model  \eqref{App_mixture}
\begin{eqnarray*}
E_{F_i}[\{u_i(X)-\G_i^* X \}^2]&=& \nu_i^* u_i^2(0)+(1-\nu_i^*) E_{G_i}[\{u_i(X)-\G_i^* X \}^2], 
\end{eqnarray*}
where $E_{G_i}$ indicates that the expectation is taken with respect to $G_i$. 
Then 
\begin{eqnarray}
 \label{App_gini.non.m21}
E_{F_i}[\{u_i(X)-\G_i^* X \}^2]&=& \nu_i^* u_i^2(0)+(1-\nu_i^*)E_{G_i}\{u_i^2(X)-2\G_i^* Xu_i (X)+\G_i^{*2} X^2 \}.  
\end{eqnarray}
With the form of $u_i(x;\bnu)$ in  \eqref{App_gini.u01}, we have 
\begin{eqnarray}
 \label{App_gini.non.u0}
u_i(0)&=&u_i(0;\bnu^*)=2(1-\nu_i^*)m_i-(1-\nu^*_i)\psi_i .  
\end{eqnarray}

Combining  \eqref{App_gini.non.var1}-- \eqref{App_gini.non.u0} gives 
\begin{eqnarray}
\nonumber
&&Var_{F_i}\{u_i(X)-\G_i X \}\\
 \label{App_gini.non.var2}&=& \nu_i^*(1-\nu_i^*)^3(2m_i-\psi_i)^2+(1-\nu_i^*)E_{G_i}\{u_i^2(X)-2\G_i^* Xu_i (X)+\G_i^{*2} X^2 \}. 
\end{eqnarray}
The fact that $\mu_i=(1-\nu_i^*) m_i$  and   \eqref{App_gini.non.var2} together imply that $\sigma_i^2$ in  \eqref{App_gini.sigma.i1} has the following form: 
\begin{eqnarray*}
    \sigma_i^2 &=&\frac{  \nu_i^*(1-\nu_i^*)^3(2m_i-\psi_i)^2+(1-\nu_i^*)E_{G_i}\{u_i^2(X)-2\G_i^* Xu_i (X)+\G_i^{*2} X^2 \} }{w_i(1-\nu_i^*)^2m_i^2}\\
    &=&\frac{  \nu_i^*(1-\nu_i^*)^3(2m_i-\psi_i)^2+(1-\nu_i^*)E_{G_i}\{u_i^2(X)-2\G_i^* Xu_i (X)+\G_i^{*2} X^2 \} }{w_i(1-\nu_i^*)^2m_i^2}\\
    &=&\frac{ E_{G_i}\{u_i^2(X)-2\G_i^* Xu_i (X)+\G_i^{*2} X^2 \} }{w_i(1-\nu_i^*) m_i^2}
    +
   \frac{  \nu_i^* (1-\G_i^*)^2 }{w_i(1-\nu_i^*) }, 
\end{eqnarray*}
where in the last step, we have used the fact that  $\G_i^*
= (2 \nu_i^* - 1) + (1-\nu_i^*) {\psi_i}/{m_i}$.

Under the DRM  \eqref{App_drm} and since $\Delta^*\rho^*=w_1(1-\nu_1^*)$ and $\Delta^*(1-\rho^*)=w_0(1-\nu_0^*)$,
we further have 
\begin{eqnarray*}
    \sigma_0^2 &=&
\frac{ E_{0}\{u_0^2(X)-2\G_0^* Xu_0 (X)+\G_0^{*2} X^2 \} }{\Delta^*(1-\rho^*)  m_0^2}
    +
   \frac{  \nu_0^* (1-\G_0^*)^2 }{\Delta^*(1-\rho^*) } 
\end{eqnarray*}
and 
\begin{eqnarray*}
    \sigma_1^2 &=&
\frac{ E_{0}[\omega(X) \{u_1^2(X)-2\G_1^* Xu_1 (X)+\G_1^{*2} X^2 \} ]}{\Delta^* \rho^*  m_1^2}
    +
   \frac{  \nu_1^* (1-\G_1^*)^2 }{\Delta^* \rho^* } .
\end{eqnarray*}

Recall that 
$$
\J_0=\left(-\frac{\G_0^*}{m_0}, \frac{1}{m_0}\right)^\top, ~~
\J_1=\left(-\frac{\G_1^*}{m_1}, \frac{1}{m_1}\right)^\top
$$
and 
$$
\tilde{\bu}_0(x)=-\rho^*\big(x,u_0(x)\big)^\top,
~~
\tilde{\bu}_1(x)=(1-\rho^*)\big(x,u_1(x)\big)^\top
.
$$
After some algebra work, we get 
\begin{eqnarray}
 \label{App_sigma_02}
\sigma_0^2 
=\frac{1}{\Delta^*(\rho^*)^2(1-\rho^*)}\J_0^\top
E_0\{\tilde{\bu}_0(X)\tilde{\bu}_0(X)^\top\}\J_0+ \frac{\nu_0^*(\G_0^*-1)^2}{\Delta^*(1-\rho^*)}
\end{eqnarray}
and  
\begin{eqnarray}
 \label{App_sigma_12}
\sigma_1^2
= \frac{1}{\Delta^*\rho^*(1-\rho^*)^2}\J_1^\top E_0\{\omega(X)\tilde{\bu}_1(X)\tilde{\bu}_1(X)^\top\} \J_1+ \frac{\nu_1^*(\G_1^*-1)^2}{\Delta^*\rho^*}.
\end{eqnarray}

Substituting   \eqref{App_sigma_02} and  \eqref{App_sigma_12} into  \eqref{App_bSigma_non} gives the asymptotic variance $\bSigma_{non}$ as 
\begin{eqnarray}
\bSigma_{non} = 
   \J\bSigma_{np1}\J^\top + \diag\lb \frac{\nu_0^*(1 - \G_0^*)^2}{\Delta^*(1-\rho^*)},\frac{\nu_1^*(1 - \G_1^*)^2}{\Delta^*\rho^*}\rb,
\end{eqnarray}
with 
\begin{equation*}
    \bSigma_{np1} = 
   \frac{1}{\Delta^*\rho^*(1-\rho)}\diag\lb\frac{E_0\{\tilde{\bu}_0(X)\tilde{\bu}_0(X)^\top\}}{\rho^*},\frac{E_0\{\omega(X)\tilde{\bu}_1(X)\tilde{\bu}_1(X)^\top\}}{1-\rho^*}\rb.
\end{equation*}
Hence, $\bSigma_{non}$ has the form claimed in (a). 

We now move to (b). 
Since $\bu(x;\bnu^*,\btheta^*) = (-(\rho^*)^{-1}\tilde{\bu}_0(x)^\top,(1-\rho^*)^{-1}\tilde{\bu}_1(x)^\top)^\top$, after some algebra, we find that
\begin{eqnarray*}
&&\frac{1}{\Delta^*}E_0\left\{\frac{\bu(X;\bnu^*,\btheta^*)\bu(X;\bnu^*,\btheta^*)^\top}{h(X)}\right\} \\
&=& \frac{1}{\Delta^*(\rho^*)^2(1-\rho^*)}\ba{cc}E_0\{h_0(X)\tilde{\bu}_0(X)\tilde{\bu}_0(X)^\top\}& -E_0\{h_1\tilde{\bu}_0(X)\tilde{\bu}_1(X)^\top\} \\
-E_0\{h_1(X)\tilde{\bu}_1(X)\tilde{\bu}_0(X)^\top\}& \frac{\rho^*}{1-\rho^*}E_0\{h_1(X)\omega(X)\tilde{\bu}_1(X)\tilde{\bu}_1(X)^\top\}\ea\\
&=&\bSigma_{np1} - \frac{1}{\Delta^*(\rho^*)^2(1-\rho)}E_0\{h_1(X)\tilde{\bu}(X)\tilde{\bu}(X)^\top\}.
\end{eqnarray*}
Together with the expression for $\bSigma$ in  \eqref{App_final_bSigma}, it follows that
\begin{equation*}
 \bSigma_{non}   - \bSigma = 
 \frac{1}{\Delta^*(\rho^*)^2(1-\rho^*)}\J \Bigg[ E_0\{h_1(X)\tilde{\bu}(X)\tilde{\bu}(X)^\top\} -\Delta^*(1-\rho^*) \B \Bigg] \J^\top,
\end{equation*}
where 
\begin{equation}
 \label{App_gini.Bmatrix}
\B = E_0\{h_1(X)\tilde{\bu}(X)\Q(X)^\top\}\A_{\btheta}^{-1} E_0\{h_1(X)\Q(X)\tilde{\bu}(X)^\top\}.
\end{equation} 

Let $ \D(a) = \left(\D_0(a)^\top,\D_1(a)^\top\right)^\top$ for $a>0$ with 
\begin{equation*}
   \D_i(a)=\tilde{\bu}_i(x)- \Delta^*(1-\rho^*) E_0\left\{h_1(X)\tilde{\bu}_i(X)\Q(X)^\top\right\} \A_{\bstheta}^{-1} \Q(a), ~ i =0,1.
\end{equation*}
Recall that 
$$
\A_{\bstheta} =
\Delta^*(1-\rho^*)E_0\left[h_1(X) \Q(X)\Q^\top(X) \right].
$$
It can be verified that for $i,j \in\{0,1\}$,
\begin{eqnarray*}
&&E_0\{h_1(X)\D_i(X)\D_j(X)^\top\} \\
&=& E_0\{h_1(X)\tilde{\bu}_i(X)\tilde{\bu}_j(X)^\top\} \\
&&- \Delta^*(1-\rho^*)E_0\{h_1(X)\tilde{\bu}_i(X)\Q(X)^\top\}\A_{\bstheta}^{-1}E_0\{h_1(X)\Q(X)
\tilde{\bu}_j(X)^\top\}.
\end{eqnarray*}

Recall that 
\begin{eqnarray*}
&\tilde{\bu}(X) = (\tilde{\bu}_0(X)^\top,\tilde{\bu}_1(X)^\top)^\top
\end{eqnarray*}
and 
$\B$ is given in  \eqref{App_gini.Bmatrix}. 
Then, 
\begin{equation*}
    \bSigma_{non}   - \bSigma = 
 \frac{1}{\Delta^*(\rho^*)^2(1-\rho^*)}\J  E_0\{h_1(X)\D(X)\D(X)^\top\}  \J^\top,
\end{equation*}
as claimed in (b). This completes the proof. 

\section{Proof of Theorem 4}
The result in Theorem 4 is a direct consequence of applying the Delta method  and the results in Theorems 1 and 2.

\section{Additional simulation results}

\subsection{Results for point estimator}
Tables  \ref{App_point_chi_supp} and  \ref{App_point_exp_supp} present the additional simulated results for the point estimators of the Gini indices $\G_0$, $\G_1$, and their difference $\G_0 - \G_1$ under different distributional settings. 
The general trends are similar to those in the main paper. 
The DRM method always gives the smallest mean square errors (MSEs). 

\begin{table}[!htt]
  \centering
  \footnotesize	
  \tabcolsep 2mm
  \caption{Bias ($\times 1000$) and MSE ($\times 1000$) for point estimators ($\chi^2$ distributions).}
   \label{App_point_chi_supp}
    \begin{tabular}{ccccccccc}
    \hline
              &       &       & \multicolumn{2}{c}{$\G_0$} & \multicolumn{2}{c}{$\G_1$} & \multicolumn{2}{c}{$\G_0-\G_1$} \\
          \cline{4-9}
    $(n_0,n_1)$ & $\bnu$   &       & Bias  & MSE   & Bias  & MSE   & Bias  & MSE \\
    \hline
    (100,100) & (0.1,0.3) & EMP   & 6.13  & 0.98  & 7.43  & 1.30  & -1.30 & 2.23 \\
          &       & JEL   & 0.96  & 0.96  & 3.08  & 1.28  & -2.13 & 2.28 \\
          &       & DRM   & 2.51  & 0.67  & 3.79  & 1.10  & -1.28 & 1.40 \\[1mm]
          & (0.6,0.4) & EMP   & 7.15  & 1.14  & 6.00  & 1.27  & 1.15  & 2.44 \\
          &       & JEL   & 4.90  & 1.14  & 2.27  & 1.27  & 2.63  & 2.49 \\
          &       & DRM   & 2.71  & 0.94  & 3.48  & 1.18  & -0.77 & 1.98 \\[2mm]
    (300,300) & (0.1,0.3) & EMP   & 1.40  & 0.31  & 2.96  & 0.41  & -1.56 & 0.72 \\
          &       & JEL   & -0.33 & 0.31  & 1.51  & 0.40  & -1.83 & 0.72 \\
          &       & DRM   & 0.75  & 0.23  & 1.38  & 0.35  & -0.63 & 0.46 \\[1mm]
          & (0.6,0.4) & EMP   & 2.63  & 0.38  & 1.58  & 0.43  & 1.05  & 0.80 \\
          &       & JEL   & 1.87  & 0.38  & 0.33  & 0.43  & 1.54  & 0.81 \\
          &       & DRM   & 1.02  & 0.32  & 0.78  & 0.41  & 0.25  & 0.66 \\
          \hline
    \end{tabular}%
\end{table}%

\begin{table}[!htt]
  \centering
  \footnotesize	
  \tabcolsep 2mm
  \caption{Bias ($\times 1000$) and MSE ($\times 1000$) for point estimators (exponential distributions).}
   \label{App_point_exp_supp}
    \begin{tabular}{ccccccccc}
    \hline
              &       &       & \multicolumn{2}{c}{$\G_0$} & \multicolumn{2}{c}{$\G_1$} & \multicolumn{2}{c}{$\G_0-\G_1$} \\
          \cline{4-9}
    $(n_0,n_1)$ & $\bnu$   &       & Bias  & MSE   & Bias  & MSE   & Bias  & MSE \\
    \hline
    (100,100) & (0.1,0.3) & EMP   & 4.83  & 0.99  & 4.12  & 1.08  & 0.70  & 1.95 \\
          &       & JEL   & 0.33  & 0.99  & 0.63  & 1.09  & -0.30 & 1.98 \\
          &       & DRM   & 1.63  & 0.88  & 1.61  & 0.77  & 0.02  & 1.22 \\[1mm]
          & (0.6,0.4) & EMP   & 6.12  & 0.95  & 3.94  & 1.12  & 2.19  & 2.04 \\
          &       & JEL   & 4.17  & 0.95  & 0.95  & 1.13  & 3.22  & 2.09 \\
          &       & DRM   & 2.34  & 0.83  & 2.09  & 0.92  & 0.25  & 1.51 \\[2mm]
    (300,300) & (0.1,0.3) & EMP   & 1.52  & 0.33  & 2.26  & 0.37  & -0.75 & 0.66 \\
          &       & JEL   & 0.02  & 0.33  & 1.10  & 0.36  & -1.08 & 0.67 \\
          &       & DRM   & 0.70  & 0.29  & 0.96  & 0.25  & -0.25 & 0.40 \\[1mm]
          & (0.6,0.4) & EMP   & 2.03  & 0.31  & 1.26  & 0.38  & 0.77  & 0.69 \\
          &       & JEL   & 1.37  & 0.31  & 0.27  & 0.38  & 1.11  & 0.70 \\
          &       & DRM   & 0.83  & 0.29  & 0.45  & 0.30  & 0.39  & 0.51 \\
          \hline
    \end{tabular}%
\end{table}%

\subsection{Results for confidence intervals}

Tables  \ref{App_CI_chi_supp} and  \ref{App_CI_exp_supp} contain the the complete results for the confidence intervals (CIs) of $\G_0$ and $\G_1$ under different distributional settings. 
NL-DRM and BL-DRM refer to the Wald-type CIs for $\G_0$ or $\G_1$ using the logit transformation under the DRM and the corresponding bootstrap-t CIs.
The additional results for the CIs of $\G_0-\G_1$ are shown in Table~ \ref{App_CI_diff_supp}. 
Again, the general patterns are similar to those in the main paper. 
The NA-DRM CIs provide accurate coverage probabilities (CPs) in all situations and have shorter average lengths (ALs)
than the existing nonparametric methods. 
Further, the bootstrap method and logit transformation do not help to improve the coverage accuracy. 
Hence, we recommend using the NA-DRM CI. 

\begin{table}[!htt]
  \centering
  \footnotesize	
  \tabcolsep 2mm
  \caption{Coverage probability (CP\%) and average length (AL) of CIs ($\chi^2$ distributions).}
   \label{App_CI_chi_supp}
    \resizebox{13cm}{10.5cm}{ \begin{tabular}{cc cccccccc}
    \hline
          & \multicolumn{1}{r}{} & \multicolumn{4}{c}{(100,100)}       & \multicolumn{4}{c}{(300,300)} \\
    \cline{3-10}
          & \multicolumn{1}{c}{} & \multicolumn{2}{c}{$\G_0$} & \multicolumn{2}{c}{$\G_1$} & \multicolumn{2}{c}{$\G_0$} & \multicolumn{2}{c}{$\G_1$} \\
    $\bnu$    & \multicolumn{1}{c}{} & \multicolumn{1}{c}{CP} & \multicolumn{1}{c}{AL} & \multicolumn{1}{c}{CP} & \multicolumn{1}{c}{AL} & \multicolumn{1}{c}{CP} & \multicolumn{1}{c}{AL} & \multicolumn{1}{c}{CP} & \multicolumn{1}{c}{AL} \\
    \hline
    (0,0) & NA-EMP & 93.85 & 0.100 & 94.20 & 0.092 & 94.60 & 0.059 & 94.80 & 0.054 \\
          & BT-EMP & 94.10 & 0.103 & 94.75 & 0.094 & 94.85 & 0.059 & 95.05 & 0.054 \\
          & EL    & 93.85 & 0.100 & 94.20 & 0.091 & 94.55 & 0.059 & 94.80 & 0.054 \\
          & BT-EL & 94.45 & 0.103 & 95.10 & 0.095 & 94.90 & 0.059 & 94.95 & 0.054 \\
          & JEL   & 94.45 & 0.102 & 94.85 & 0.094 & 94.70 & 0.059 & 95.15 & 0.054 \\
          & AJEL  & 94.80 & 0.105 & 95.50 & 0.096 & 94.90 & 0.060 & 95.30 & 0.055 \\
          & NA-DRM & 95.25 & 0.074 & 94.65 & 0.078 & 94.70 & 0.043 & 94.70 & 0.045 \\
          & BT-DRM & 95.55 & 0.075 & 95.00 & 0.079 & 94.55 & 0.043 & 94.55 & 0.046 \\
          & NL-DRM & 95.35 & 0.074 & 94.50 & 0.077 & 94.75 & 0.043 & 94.75 & 0.045 \\
          & BL-DRM & 95.45 & 0.075 & 94.80 & 0.079 & 94.50 & 0.043 & 94.55 & 0.045 \\[1.5mm]
    (0.1,0.3) & NA-EMP & 94.00 & 0.116 & 93.85 & 0.134 & 95.00 & 0.068 & 95.10 & 0.079 \\
          & BT-EMP & 94.80 & 0.119 & 95.05 & 0.137 & 95.35 & 0.068 & 95.10 & 0.079 \\
          & EL    & 93.90 & 0.116 & 93.95 & 0.133 & 95.00 & 0.068 & 95.05 & 0.078 \\
          & BT-EL & 95.25 & 0.119 & 94.65 & 0.139 & 95.20 & 0.068 & 95.10 & 0.080 \\
          & JEL   & 94.70 & 0.120 & 94.00 & 0.140 & 95.25 & 0.069 & 94.65 & 0.080 \\
          & AJEL  & 95.25 & 0.123 & 94.60 & 0.144 & 95.40 & 0.069 & 94.80 & 0.081 \\
          & NA-DRM & 93.60 & 0.099 & 94.80 & 0.128 & 94.60 & 0.058 & 95.25 & 0.075 \\
          & BT-DRM & 93.95 & 0.100 & 95.25 & 0.129 & 94.55 & 0.058 & 94.95 & 0.074 \\
          & NL-DRM & 93.85 & 0.099 & 95.00 & 0.128 & 94.65 & 0.058 & 95.20 & 0.075 \\
          & BL-DRM & 93.65 & 0.099 & 94.95 & 0.127 & 94.55 & 0.058 & 94.85 & 0.074 \\[1.5mm]
    (0.3,0.3) & NA-EMP & 93.80 & 0.132 & 93.65 & 0.134 & 94.60 & 0.077 & 94.05 & 0.079 \\
          & BT-EMP & 95.30 & 0.135 & 94.55 & 0.137 & 95.20 & 0.077 & 94.40 & 0.079 \\
          & EL    & 93.75 & 0.131 & 93.65 & 0.134 & 94.60 & 0.077 & 94.00 & 0.078 \\
          & BT-EL & 94.50 & 0.136 & 94.85 & 0.139 & 94.65 & 0.078 & 94.55 & 0.079 \\
          & JEL   & 94.45 & 0.137 & 93.80 & 0.141 & 94.50 & 0.078 & 94.55 & 0.080 \\
          & AJEL  & 95.35 & 0.141 & 94.20 & 0.144 & 94.80 & 0.079 & 94.80 & 0.081 \\
          & NA-DRM & 95.10 & 0.120 & 94.35 & 0.130 & 95.45 & 0.070 & 94.90 & 0.076 \\
          & BT-DRM & 95.75 & 0.121 & 94.65 & 0.130 & 95.25 & 0.070 & 94.65 & 0.075 \\
          & NL-DRM & 95.60 & 0.120 & 94.70 & 0.129 & 95.50 & 0.070 & 95.00 & 0.076 \\
          & BL-DRM & 95.30 & 0.119 & 94.60 & 0.128 & 95.10 & 0.069 & 94.60 & 0.075 \\  [1.5mm]  
    (0.6,0.4) & NA-EMP & 93.45 & 0.124 & 94.10 & 0.138 & 94.30 & 0.073 & 95.10 & 0.080 \\
          & BT-EMP & 95.85 & 0.131 & 95.00 & 0.142 & 95.55 & 0.074 & 95.10 & 0.081 \\
          & EL    & 94.00 & 0.123 & 94.10 & 0.137 & 94.45 & 0.073 & 95.05 & 0.080 \\
          & BT-EL & 95.35 & 0.130 & 94.95 & 0.143 & 94.90 & 0.075 & 95.30 & 0.082 \\
          & JEL   & 92.90 & 0.133 & 94.15 & 0.145 & 93.35 & 0.075 & 94.90 & 0.082 \\
          & AJEL  & 93.40 & 0.137 & 94.90 & 0.149 & 93.60 & 0.075 & 95.05 & 0.083 \\
          & NA-DRM & 94.60 & 0.119 & 95.05 & 0.137 & 95.30 & 0.069 & 95.15 & 0.080 \\
          & BT-DRM & 94.95 & 0.120 & 95.45 & 0.137 & 95.55 & 0.069 & 94.95 & 0.078 \\
          & NL-DRM & 95.00 & 0.119 & 95.25 & 0.136 & 95.60 & 0.069 & 95.15 & 0.080 \\
          & BL-DRM & 94.45 & 0.116 & 94.85 & 0.134 & 95.10 & 0.068 & 94.80 & 0.078 \\ [1.5mm]
     (0.7,0.7) & NA-EMP & 92.20 & 0.113 & 92.95 & 0.119 & 94.90 & 0.067 & 93.90 & 0.070 \\
          & BT-EMP & 96.75 & 0.122 & 96.55 & 0.128 & 96.30 & 0.068 & 95.40 & 0.072 \\
          & EL    & 92.35 & 0.111 & 92.90 & 0.117 & 95.15 & 0.067 & 93.75 & 0.070 \\
          & BT-EL & 94.70 & 0.120 & 95.30 & 0.127 & 95.75 & 0.069 & 94.55 & 0.072 \\
          & JEL   & 90.75 & 0.123 & 90.80 & 0.129 & 94.00 & 0.069 & 93.00 & 0.072 \\
          & AJEL  & 91.35 & 0.127 & 91.55 & 0.133 & 94.25 & 0.070 & 93.10 & 0.073 \\
          & NA-DRM & 94.50 & 0.111 & 94.85 & 0.121 & 95.10 & 0.065 & 95.20 & 0.071 \\
          & BT-DRM & 95.40 & 0.113 & 95.90 & 0.123 & 95.45 & 0.064 & 95.60 & 0.070 \\
          & NL-DRM & 95.70 & 0.111 & 96.05 & 0.121 & 96.05 & 0.065 & 96.10 & 0.071 \\
          & BL-DRM & 94.60 & 0.108 & 95.05 & 0.118 & 95.30 & 0.063 & 95.20 & 0.069 \\
            \hline
    \end{tabular}}%
\end{table}%

\begin{table}[!htt]
  \centering
  \footnotesize	
  \tabcolsep 2mm
  \caption{Coverage probability (CP\%) and average length (AL) of CIs (exponential distributions).}
   \label{App_CI_exp_supp}
  \resizebox{13cm}{10.5cm}{  \begin{tabular}{cc cccccccc}
    \hline
          & \multicolumn{1}{r}{} & \multicolumn{4}{c}{(100,100)}       & \multicolumn{4}{c}{(300,300)} \\
    \cline{3-10}
          & \multicolumn{1}{c}{} & \multicolumn{2}{c}{$\G_0$} & \multicolumn{2}{c}{$\G_1$} & \multicolumn{2}{c}{$\G_0$} & \multicolumn{2}{c}{$\G_1$} \\
    $\bnu$    & \multicolumn{1}{c}{} & \multicolumn{1}{c}{CP} & \multicolumn{1}{c}{AL} & \multicolumn{1}{c}{CP} & \multicolumn{1}{c}{AL} & \multicolumn{1}{c}{CP} & \multicolumn{1}{c}{AL} & \multicolumn{1}{c}{CP} & \multicolumn{1}{c}{AL} \\
    \hline
    (0,0) & NA-EMP & 93.85 & 0.110 & 93.50 & 0.111 & 94.65 & 0.065 & 94.45 & 0.065 \\
          & BT-EMP & 94.35 & 0.115 & 94.05 & 0.115 & 94.75 & 0.065 & 94.75 & 0.065 \\
          & EL    & 93.90 & 0.110 & 93.50 & 0.110 & 94.65 & 0.065 & 94.55 & 0.065 \\
          & BT-EL & 94.50 & 0.113 & 94.00 & 0.113 & 94.80 & 0.065 & 94.60 & 0.065 \\
          & JEL   & 94.35 & 0.113 & 93.90 & 0.113 & 94.90 & 0.065 & 94.55 & 0.065 \\
          & AJEL  & 94.95 & 0.115 & 94.35 & 0.116 & 95.10 & 0.066 & 94.75 & 0.066 \\
          & NA-DRM & 94.80 & 0.100 & 94.05 & 0.079 & 93.95 & 0.059 & 95.20 & 0.045 \\
          & BT-DRM & 94.45 & 0.104 & 94.75 & 0.079 & 93.65 & 0.060 & 94.95 & 0.045 \\
          & NL-DRM & 94.90 & 0.100 & 94.10 & 0.078 & 93.95 & 0.059 & 95.25 & 0.045 \\
          & BL-DRM & 94.20 & 0.103 & 94.50 & 0.079 & 93.55 & 0.059 & 94.85 & 0.045 \\[1.5mm]
    (0.1,0.3) & NA-EMP & 93.45 & 0.119 & 94.30 & 0.128 & 95.25 & 0.070 & 94.50 & 0.075 \\
          & BT-EMP & 94.25 & 0.123 & 95.50 & 0.132 & 95.00 & 0.071 & 94.80 & 0.075 \\
          & EL    & 93.45 & 0.119 & 94.45 & 0.127 & 95.25 & 0.070 & 94.40 & 0.075 \\
          & BT-EL & 94.40 & 0.122 & 95.20 & 0.131 & 95.35 & 0.071 & 94.60 & 0.076 \\
          & JEL   & 94.30 & 0.122 & 94.85 & 0.133 & 95.35 & 0.071 & 94.75 & 0.076 \\
          & AJEL  & 94.80 & 0.126 & 95.20 & 0.136 & 95.40 & 0.071 & 95.05 & 0.076 \\
          & NA-DRM & 94.70 & 0.114 & 94.75 & 0.109 & 95.35 & 0.067 & 94.95 & 0.063 \\
          & BT-DRM & 94.35 & 0.116 & 95.15 & 0.109 & 95.10 & 0.066 & 94.85 & 0.063 \\
          & NL-DRM & 94.90 & 0.113 & 95.10 & 0.108 & 95.45 & 0.067 & 95.10 & 0.063 \\
          & BL-DRM & 93.90 & 0.114 & 94.70 & 0.108 & 94.85 & 0.066 & 94.85 & 0.062 \\ [1.5mm]
    (0.3,0.3) & NA-EMP & 93.55 & 0.127 & 94.25 & 0.128 & 94.55 & 0.075 & 93.45 & 0.075 \\
          & BT-EMP & 95.35 & 0.132 & 94.90 & 0.132 & 94.70 & 0.075 & 93.80 & 0.075 \\
          & EL    & 93.60 & 0.126 & 94.10 & 0.127 & 94.70 & 0.075 & 93.30 & 0.075 \\
          & BT-EL & 94.75 & 0.131 & 94.85 & 0.132 & 95.05 & 0.076 & 93.70 & 0.076 \\
          & JEL   & 93.80 & 0.132 & 94.55 & 0.133 & 94.65 & 0.076 & 93.65 & 0.076 \\
          & AJEL  & 94.50 & 0.136 & 95.15 & 0.136 & 95.00 & 0.076 & 93.80 & 0.076 \\
          & NA-DRM & 95.55 & 0.124 & 94.95 & 0.112 & 95.15 & 0.073 & 94.60 & 0.065 \\
          & BT-DRM & 95.45 & 0.125 & 95.30 & 0.112 & 94.60 & 0.072 & 94.60 & 0.064 \\
          & NL-DRM & 95.65 & 0.123 & 95.00 & 0.111 & 95.25 & 0.073 & 94.90 & 0.064 \\
          & BL-DRM & 95.25 & 0.122 & 95.15 & 0.110 & 94.45 & 0.071 & 94.55 & 0.064 \\[1.5mm]
    (0.6,0.4) & NA-EMP & 92.70 & 0.115 & 94.05 & 0.127 & 93.65 & 0.068 & 94.30 & 0.075 \\
          & BT-EMP & 95.50 & 0.123 & 95.35 & 0.132 & 94.90 & 0.070 & 94.80 & 0.075 \\
          & EL    & 93.10 & 0.113 & 94.05 & 0.126 & 93.75 & 0.068 & 94.40 & 0.074 \\
          & BT-EL & 94.30 & 0.121 & 95.20 & 0.132 & 94.15 & 0.070 & 94.75 & 0.075 \\
          & JEL   & 92.45 & 0.124 & 94.55 & 0.133 & 93.45 & 0.070 & 94.50 & 0.076 \\
          & AJEL  & 92.75 & 0.127 & 95.15 & 0.137 & 93.85 & 0.070 & 94.50 & 0.076 \\
          & NA-DRM & 94.95 & 0.116 & 95.30 & 0.118 & 94.70 & 0.068 & 94.60 & 0.068 \\
          & BT-DRM & 95.45 & 0.116 & 95.50 & 0.118 & 94.30 & 0.066 & 94.50 & 0.067 \\
          & NL-DRM & 96.05 & 0.116 & 95.80 & 0.118 & 95.10 & 0.068 & 94.85 & 0.068 \\
          & BL-DRM & 94.75 & 0.112 & 95.05 & 0.116 & 94.05 & 0.066 & 94.35 & 0.067 \\ [1.5mm]
    (0.7,0.7) & NA-EMP & 91.40 & 0.104 & 92.15 & 0.105 & 94.60 & 0.062 & 94.55 & 0.062 \\
          & BT-EMP & 96.30 & 0.114 & 95.70 & 0.115 & 95.85 & 0.064 & 95.50 & 0.064 \\
          & EL    & 92.05 & 0.102 & 92.35 & 0.102 & 94.70 & 0.062 & 94.55 & 0.062 \\
          & BT-EL & 95.00 & 0.110 & 94.20 & 0.111 & 95.40 & 0.064 & 95.25 & 0.064 \\
          & JEL   & 90.40 & 0.114 & 91.00 & 0.114 & 94.30 & 0.064 & 93.95 & 0.064 \\
          & AJEL  & 90.85 & 0.117 & 91.60 & 0.118 & 94.40 & 0.064 & 94.05 & 0.064 \\
          & NA-DRM & 94.65 & 0.109 & 93.90 & 0.101 & 96.10 & 0.064 & 95.40 & 0.059 \\
          & BT-DRM & 95.80 & 0.109 & 95.55 & 0.102 & 95.65 & 0.062 & 95.95 & 0.058 \\
          & NL-DRM & 97.05 & 0.109 & 95.75 & 0.102 & 96.75 & 0.065 & 95.95 & 0.059 \\
         & BL-DRM & 94.40 & 0.104 & 94.50 & 0.098 & 95.35 & 0.061 & 95.50 & 0.057 \\
          \hline
    \end{tabular}}%
\end{table}%

\begin{table}[!htt]
  \centering
  \footnotesize	
  \tabcolsep 2mm
  \caption{Coverage probability (CP\%) and average length (AL) of CIs for the difference $\G_0 - \G_1$.}
   \label{App_CI_diff_supp}
    \begin{tabular}{cc cccccccc}
    \hline
          & \multicolumn{1}{r}{} & \multicolumn{4}{c}{(100,100)}     & \multicolumn{4}{c}{(300,300)} \\
          \cline{3-10}
          & \multicolumn{1}{c}{} & \multicolumn{2}{c}{$\chi^2$} & \multicolumn{2}{c}{$Exp$} & \multicolumn{2}{c}{$\chi^2$} & \multicolumn{2}{c}{$Exp$} \\
          & \multicolumn{1}{c}{} & \multicolumn{1}{c}{CP} & \multicolumn{1}{c}{AL} & \multicolumn{1}{c}{CP} & \multicolumn{1}{c}{AL} & \multicolumn{1}{c}{CP} & \multicolumn{1}{c}{AL} & \multicolumn{1}{c}{CP} & \multicolumn{1}{c}{AL} \\
    \hline
    (0.1,0.3) & NA-EMP & 94.45 & 0.178 & 94.55 & 0.104 & 95.20 & 0.175 & 94.69 & 0.103 \\
          & BT-EMP & 95.10 & 0.181 & 94.55 & 0.104 & 95.15 & 0.179 & 94.74 & 0.103 \\
          & JEL   & 95.45 & 0.190 & 95.15 & 0.106 & 96.40 & 0.187 & 95.15 & 0.105 \\
          & AJEL  & 95.60 & 0.193 & 95.25 & 0.107 & 96.60 & 0.190 & 95.25 & 0.106 \\
          & NA-DRM & 94.40 & 0.146 & 95.65 & 0.085 & 94.85 & 0.138 & 94.89 & 0.080 \\
          & BT-DRM & 93.45 & 0.143 & 95.10 & 0.083 & 93.80 & 0.135 & 94.34 & 0.080 \\ [1.5mm]
    (0.6,0.4) & NA-EMP & 94.60 & 0.186 & 94.30 & 0.109 & 93.84 & 0.172 & 95.00 & 0.101 \\
          & BT-EMP & 94.90 & 0.190 & 94.70 & 0.109 & 94.39 & 0.177 & 95.10 & 0.102 \\
          & JEL   & 96.25 & 0.205 & 94.90 & 0.112 & 96.85 & 0.192 & 96.30 & 0.105 \\
          & AJEL  & 96.40 & 0.208 & 95.00 & 0.113 & 97.10 & 0.196 & 96.40 & 0.105 \\
          & NA-DRM & 95.35 & 0.175 & 95.25 & 0.102 & 95.40 & 0.150 & 95.10 & 0.088 \\
          & BT-DRM & 93.95 & 0.168 & 94.05 & 0.098 & 94.59 & 0.147 & 94.75 & 0.086 \\
    \hline
    \end{tabular}%
\end{table}%

\clearpage 

\end{document}